\documentclass[twocolumn]{autart}
\usepackage{amsmath,amsfonts,amssymb,amscd,mathrsfs,array,subeqnarray}
\allowdisplaybreaks
\usepackage{lipsum}
\begin{document}
\begin{frontmatter}

\title{Two person non-zero-sum linear-quadratic differential game with Markovian jumps in infinite horizon\thanksref{footnoteinfo}}
\thanks[footnoteinfo]{This paper was not presented at any conference. Corresponding author Xin Zhang.} 

\author[a]{Fan Wu}\ead{wfyy121107@163.com},
\author[b]{Xun Li}\ead{li.xun@polyu.edu.hk},
\author[a]{Xin Zhang}\ead{x.zhang.seu@gmail.com}

\address[a]{School of Mathematics, Southeast University, Nanjing 211189, China}  
\address[b]{Department of Applied Mathematics, Hong Kong Polytechnic University, Hong Kong, China}             



\begin{abstract}
  This paper investigates an inhomogeneous non-zero-sum linear-quadratic (LQ, for short) differential game problem whose state process and cost functional are regulated by a Markov chain. Under the $L^2$ stabilizability framework, we first provide a sufficient condition to ensure the $L^2$-integrability of the state process and study a class of linear backward stochastic differential equation (BSDE, for short) in infinite horizon. Then, we seriously discuss the LQ  problem and show that the closed-loop optimal control is characterized by the solutions to coupled algebra Riccati equations (CAREs, for short) with some stabilizing conditions and a linear BSDE. Based on those results, we further analyze the non-zero-sum stochastic differential game problem and give the closed-loop Nash equilibrium through the solution to a system of two cross-coupled CAREs and two cross-coupled BSDEs. Finally, some related numerical examples are presented to illustrate the results obtained in the paper.
\end{abstract}
\end{frontmatter}

\section{Introduction}\label{section-introduction}
Let $(\Omega,\mathcal{F},\mathbb{F},\mathbb{P} )$ be a complete filtered probability space on which a continuous time, finite-state, irreducible Markov chain $\alpha$ are defined with $\mathbb{F}=\{\mathcal{F}_{t}\}_{t\geq 0}$ being its natural filtration augmented by all $\mathbb{P}$-null sets in $\mathcal{F}$. The state space of the Markov chain $\alpha$ is denoted by $\mathcal{S}:=\left\{1,2,...,D\right\}$, where $D$ is a finite natural number. We define $\Theta(\alpha_{t})\triangleq\sum_{i=1}^{D}\Theta(i)\mathbb{I}_{(\alpha_{t}=i)}(t)$ for any given D-dimension vector $\mathbf{\Theta}=[\Theta(1),\Theta(2),...,\Theta(D)]$, where $\mathbb{I}_{A}$ is the indicator function. For any Euclidean space $\mathbb{H}$,
let $L_{\mathbb{F}}^{2,loc}(\mathbb{H})$ be the set of $\mathbb{H}$-valued, $\mathbb{F}$-progressively measurable process with $\mathbb{E} \int_{0}^{T}|\varphi(s)|^{2} ds<\infty$ for all $T>0$ and let $L_{\mathbb{F}}^{2}(\mathbb{H})$ ($L_{\mathcal{P}}^{2}(\mathbb{H})$)  be the set of $\mathbb{H}$-valued, $\mathbb{F}$-progressively measurable ($\mathbb{F}$-predictable) process with $\mathbb{E} \int_{0}^{\infty}|\varphi(s)|^{2} ds<\infty$. Let $\mathbb{S}^n$ ($\mathbb{S}_{+}^n$, $\overline{\mathbb{S}_{+}^n}$) represents the set of all $n \times n$ symmetric matrices (positive-definite matrices, positive-semidefinite matrices). For $M, N \in \mathbb{S}^n$, we write $M \geq N$ (respectively, $M>N$) if $M-N$ is positive-semidefinite (respectively, positive-definite).

We consider a two-person non-zero-sum LQ differential game problem with Markovian jumps. The controlled system is given by
 \begin{equation}\label{GLQ-state}
   \left\{
   \begin{aligned}
   \dot{X}(t)&\!=\!A(\alpha_{t})X(t)\!+\!B_{1}(\alpha_{t})u_{1}(t)\!+\!B_{2}(\alpha_{t})u_{2}(t)\!+\!b(t),\\
   X(0)&=x,\quad \alpha(0)=i,
   \end{aligned}
   \right.
 \end{equation}
 and the cost functional of two players are defined as:
\begin{equation}\label{GLQ-cost}
\begin{aligned}
    &\quad J_{k}\left(x,i;u_{1},u_{2}\right)\\
    &\triangleq \mathbb{E}\int_{0}^{\infty}\!\left[\!
    \left<\!
    \left(\!
    \begin{matrix}
    Q^{k}(\alpha)\!&S_{1}^{k}(\alpha)^{\top}\!&S_{2}^{k}(\alpha)^{\top}\\
    S_{1}^{k}(\alpha)\!&R_{11}^{k}(\alpha)\!&R_{12}^{k}(\alpha)\\
    S_{2}^{k}(\alpha)\!&R_{21}^{k}(\alpha)\!&R_{22}^{k}(\alpha)
    \end{matrix}
    \!\right)\!
    \left(\!
    \begin{matrix}
    X\\
    u_{1}\\
    u_{2}
    \end{matrix}
    \!\right)\!,
    \left(\!
    \begin{matrix}
    X\\
    u_{1}\\
    u_{2}
    \end{matrix}
    \!\right)\!
    \right>\!\right.\\
    &\quad\left.+2\left<\!\left(\!\begin{matrix}
    q^{k}\\
    \rho_{1}^{k}\\
    \rho_{2}^{k}
    \end{matrix}
    \!\right)\!,
    \left(\!\begin{matrix}
    X\\
    u_{1}\\
    u_{2}
    \end{matrix}
    \!\right)\!\right>\!\right]dt,\quad k=1,2.
  \end{aligned}
\end{equation}
In the above, we suppress the variable $t$ on the right-hand side of \eqref{GLQ-cost} for convenience, and the suppressed notations will be frequently used in the rest of this paper when no confusion arises.
$X\equiv X(\cdot;x,i,u_{1},u_{2})$, valued in $\mathbb{R}^{n}$, is called the \emph{state process} and  $u_{k}$, valued in $\mathbb{R}^{m_{k}}$, is called the \emph{control process} of Player $k$. we assume that the coefficients in the state process and cost functional satisfy:
\begin{align*}
&b\in L_{\mathbb{F}}^{2}(\mathbb{R}^{n}),\quad q^{k}\in L_{\mathbb{F}}^{2}(\mathbb{R}^{n}),\\
&\rho_{1}^{k}\in L_{\mathbb{F}}^{2}(\mathbb{R}^{m_{1}}),\quad \rho_{2}^{k}\in L_{\mathbb{F}}^{2}(\mathbb{R}^{m_{2}}),\\
&A(i)\in\mathbb{R}^{n\times n}, \quad B_{k}(i)\in\mathbb{R}^{n\times m_{k}},\\
&Q^{k}(i) \in \mathbb{S}^{n}, \quad R_{ll}^{k}(i) \in \mathbb{S}_{+}^{m_{l}},\\
& S_{l}^{k}(i)\in\mathbb{R}^{m_{l} \times n},\quad R_{12}^{k}(i)=R_{21}^{k}(i)^{\top}\in\mathbb{R}^{m_{1} \times m_{2}},\\
&Q^{k}(i)-S_{k}^{k}(i)^{\top}R_{kk}^{k}(i)S_{k}^{k}(i)\in \overline{\mathbb{S}_{+}^n},\quad k,l\in\{1,2\}.
 \end{align*}

Clearly, for any initial state $(x,i)\in\mathbb{R}^{n}\times\mathcal{S}$ and control pair $(u_{1},u_{2})\in L_{\mathcal{P}}^{2}(\mathbb{R}^{m_{1}})\times L_{\mathcal{P}}^{2}(\mathbb{R}^{m_{2}})$, the state equation \eqref{GLQ-state} admits a unique solution $X(\cdot;x,i,u_{1},u_{2})\in L_{\mathbb{F}}^{2,loc}(\mathbb{R}^{n})$. To ensure the cost  functional \eqref{GLQ-cost} is well-defined, for any given initial value $(x,i)\in \mathbb{R}^{n}\times\mathcal{S}$, we introduce the following admissible control set:
\begin{align*}
\mathcal{U}_{ad}(x,i)\triangleq\big\{&(u_{1},u_{2})\in L_{\mathcal{P}}^{2}(\mathbb{R}^{m_{1}})\times L_{\mathcal{P}}^{2}(\mathbb{R}^{m_{2}}) \mid \\
& X(\cdot;x,i,u_{1},u_{2})\in L_{\mathbb{F}}^{2}(\mathbb{R}^{n})\big\}.
\end{align*}
A pair $(u_{1},u_{2})\in \mathcal{U}_{ad}(x,i)$ is called an admissible control pair for the initial state $(x,i)$. Then, the LQ differential game problem can be defined as follows. \\
\textbf{Problem (M-GLQ).} For any given $(x,i)\in \mathbb{R}^{n}\times \mathcal{S}$, find a $(u_{1}^{*},u_{2}^{*})\in\mathcal{U}_{ad}(x,i)$ such that
\begin{equation}\label{GLQ-value-function}
    \left\{
    \begin{aligned}
& J_{1}\left(x,i;u_{1}^{*},u_{2}^{*}\right)=\inf_{(u_{1},u_{2}^{*})\in \mathcal{U}_{ad}(x,i)}J_{1}\left(x,i;u_{1},u_{2}^{*}\right),\\
& J_{2}\left(x,i;u_{1}^{*},u_{2}^{*}\right)=\inf_{(u_{1}^{*},u_{2})\in \mathcal{U}_{ad}(x,i)}J_{2}\left(x,i;u_{1}^{*},u_{2}\right).
\end{aligned}
    \right.
\end{equation}
Any $(u_{1}^{*},u_{2}^{*})\in\mathcal{U}_{ad}(x,i)$ satisfying above condition is called an open-loop Nash equilibrium point of Problem (M-GLQ) for the initial value $(x,i)$, and the corresponding $X^{*}\equiv X(\cdot;x,i,u_{1}^{*},u_{2}^{*})$ is called an equilibrium state process. In addition, if $b$, $q^{k}$, $\rho_{1}^{k}$, $\rho_{2}^{k}=0$,  then the corresponding problem, cost functional are denoted by Problem (M-GLQ)$^0$, $J_{k}^{0}(x,i;u_{1},u_{2})$, respectively.

Modeling a dynamic system by a Markov chain can better characterize the instantaneous changes of system state and has been widely applied in various fields such as engineering, financial management, and economics; see, for examples, \cite{Wen_2023,Li-zhou-2002-ID-MLQ-F,Ji-Chizeck-1990-D-MLQ-I/F,sun_risk-sensitive_2018,Zhang-Siu-Meng-2010} and the references therein. At the same time, the differential games theory has played an essential role in the economy, finance, reinsurance, and elsewhere. Some early works on the mathematics theory of differential games include \cite{Basar-1982-AcademicPress-Gamebook,Yeung-2006-Spring-Gamebook}. Under the non-Markovian framework, Karatzas and Zamfirescu \cite{karatzas-2008-AP-Game} introduced a martingale approach to study continuous-time stochastic differential games of control and stopping. Elliott and Davis \cite{Elliott-1981-SIAM-Game} studied the two-person zero-sum Stackelberg differential game and obtained the feedback strategy of the game. Fleming and Souganidis \cite{Fleming-1989-IDMJ-Game} investigated the existence of a value for a zero-sum stochastic differential game using the dynamic programming principle method and viscosity technique. Tang and Hou \cite{Tang-2007-SIAM-Game} generalized the results in \cite{Fleming-1989-IDMJ-Game} by considering a general stochastic differential system and formulated the corresponding switching game. Regarding applying game theory in reinsurance, we refer the interested readers to \cite{chen_stochastic_2019,bai_hybrid_2022}.

It is worth mentioning that most of the literature discussed above only considers the differential game problem in finite time horizon. There are few studies considering the differential games in infinite horizon. Song et al. \cite{Song-2008-IEEE} studied a zero-sum differential game of regime-switching diffusion in which the terminal is a stopping time. They developed a numerical method using Markov chain approximation techniques and proved the existence of saddle points for the stochastic differential game. Then, Zhu et al. \cite{Zhu-Zhang-Bin-2014} further investigated an infinite horizon homogeneous LQ stochastic Nash differential games with Markovian jumps. Although they considered the diffusion model, control did not enter the diffusion term. Based on the existing results in \cite{Li-Zhou-Rami-2003-ID-MLQ-IF} as well as the technique of completing the square, they consecutively obtained the optimal control and equilibrium point with feedback representations for LQ problem and non-zero-sum Nash games.

However, both \cite{Zhu-Zhang-Bin-2014} and \cite{Li-Zhou-Rami-2003-ID-MLQ-IF} are studying infinite horizon LQ problems under the mean-square stable sense, whose initial work can be traced back to \cite{Rami-Zhou-Moore-2000-ID-LQ-IF,Rami-Zhou-2000-LMI-RE-IDLQIF}. Under such a framework, it is difficult to discuss the inhomogeneous LQ control problems in infinite horizon since constructing the closed-form optimal strategy requires us to solve a linear BSDE with infinite horizon, whose solvability is hard to obtain under the mean-square framework. Recently, Huang et al. \cite{Jianhui-Huang-2015} formulated an infinite horizon LQ problem with mean-field under the $L^{2}$-stabilizability sense. Based on this framework, Sun et al. \cite{Sun.JR_2016_IZSLQI} further studied an inhomogeneous zero-sum LQ stochastic differential game in infinite horizon. To construct the closed-loop saddle point, they studied a class of linear BSDE in infinite horizon and obtained its solvability under the $L^{2}$-stable condition. Then, based on those results, they characterized the closed-loop saddle points by the solutions to an algebraic Riccati equation with certain stabilizing conditions and a linear BSDE in infinite horizon.

In this paper, we consider a non-zero-sum inhomogeneous LQ differential game with Markovian jumps in infinite horizon. Although it seems that the only difference between this paper and \cite{Zhu-Zhang-Bin-2014}  is the addition of inhomogeneous terms, the two papers are actually quite different.
Firstly, we should formulate our problem under the $L^{2}$-stabilizability sense rather than the mean-square stable sense since the existence of an inhomogeneous term. To do this, we need to extend the framework of $L^{2}$-stabilizability control system introduced in \cite{Jianhui-Huang-2015} to the one regulated by a Markov chain. Secondly, to construct the closed-loop equilibrium point, we also need to investigate the solvability of a class of linear BSDE driven by Markov chains under the $L^{2}$-stable framework. All of those make our work significantly different from \cite{Zhu-Zhang-Bin-2014} and in turn make our work meaningful.

The rest of the paper is organized as follows. Section \ref{section-Preliminary} aims to investigate the $L^{2}$-stability of the linear system with Markovian jumps and study the solvability of a class of linear BSDE driven by a Markov chain in infinite horizon. In Section \ref{section-LQ}, we analyze an inhomogeneous LQ optimal control problem with Markovian jumps and obtain the corresponding closed-loop optimal control.  Based on those results, the non-zero-sum inhomogeneous LQ differential game is studied in Section \ref{section-GLQ}. Section \ref{section-Examples} concludes the paper by presenting two concrete examples to illustrate the results developed in the previous sections.

\section{Preliminary results} \label{section-Preliminary}
We begin by introducing some notations. The generator of Markov chain $\alpha$ is given by $\mathbf{\Pi}=\left\{\pi_{ij}\mid i,j\in\mathcal{S}\right\}$,
where $\pi_{ij}\geq 0$ for $i\neq j $ and $\pi_{ii}=-\sum_{j\neq i}\pi_{ij}$.
Let $N_{j}(t)$ be the number of jumps to state $j$ up to time $t$. Then
$$
\widetilde{N}_{j}(t)\triangleq N_{j}(t)-\int_{0}^{t}\sum_{i\neq j}\left[\pi_{ij}\mathbb{I}_{\{\alpha_{s-}=i\}}(s)\right]ds,
$$
is an $\left(\mathbb{F},\mathbb{P}\right)$-martingale for each $j\in\mathcal{S}$.
For a D-dimension vector process $\mathbf{\Gamma}=\left[\Gamma_1,\Gamma_2,...,\Gamma_D\right]$, we define
$$
 \mathbf{\Gamma}(s)\cdot d\mathbf{\widetilde{N}}(s)\triangleq\sum_{k=1}^{D}\Gamma_{k}(s)d\widetilde{N}_{k}(s).
$$
In addition, for a  Banach space $\mathbb{B}$, we denote
$$
\mathcal{D}(\mathbb{B})=\left\{\mathbf{\Lambda}=\left[\Lambda(1),\Lambda(2),...,\Lambda(D)\right] \mid \Lambda_i \in \mathbb{B}\text{, } \forall i\in \mathcal{S}\right\}.
$$
Specially, if $\mathbf{\Lambda}\in \mathcal{D}\left(\mathbb{R}^{n\times n}\right)$, then we say $\lambda$ ($\mu$) is the smallest (largest) eigenvalue of $\mathbf{\Lambda}$ (or process $\Lambda(\alpha)$) if it satisfies
$$
\lambda=\min\{\lambda_{1},\lambda_{2},...,\lambda_{D}\}\quad
\left(\mu=\max\{\mu_{1},\mu_{2},...,\mu_{D}\}\right),
$$
where $\lambda_{i}$ ($\mu_{i}$) is the smallest (largest) eigenvalue of $\Lambda(i)$,
$i\in\mathcal{S}$.

Now, we consider the following uncontrolled system:
\begin{equation}\label{Ab}
\left\{
\begin{aligned}
&\dot{X}(t)=A(\alpha_{t})X(t)+b(t), \quad t\geq 0,\\
&X(0)=x,\quad \alpha_{0}=i,
\end{aligned}
\right.
\end{equation}
and introduce the following condition:\\
\textbf{Condition (A). }
\begin{align*}
   \mathfrak{P}\triangleq \big\{&\mathbf{P}\in \mathcal{D}\left(\mathbb{S}_{+}^{n}\right)\big|P(i)A(i)+A(i)^{\top}P(i) \\
   &+\sum_{j=1}^{D}\pi_{ij}P(j)<0,\quad i\in\mathcal{S}\big\}\neq\emptyset.
\end{align*}

The following result shows that the differential equation \eqref{Ab} admits an L$^2$-integrable solution under the condition (A).
\begin{prop}\label{prop-FDE}
  Suppose Condition (A) holds and $b\in L_{\mathbb{F}}^{2}(\mathbb{R}^{n})$. Then the system \eqref{Ab} admits a unique solution $X\in L_{\mathbb{F}}^{2}(\mathbb{R}^{n})$.
\end{prop}
\textbf{Proof.}
The uniqueness is obvious. We only need to prove that the unique solution $X$ is in $L_{\mathbb{F}}^{2}(\mathbb{R}^{n})$.
 Let $\mathbf{P}\in \mathfrak{P}$ and  $\Lambda(i)\triangleq-\big[P(i)A(i)+A(i)^{\top}P(i)+\sum_{j=1}^{D}\pi_{ij}P(j)\big]>0$.
Applying It\^o's rule to $\left<P(\alpha_{t})X(t),X(t)\right>$, we have
$$
\begin{aligned}
&\mathbb{E}\left[\left<P(\alpha_{t})X(t),X(t)\right>-\left<P(i)x,x\right>\right]\\
=&\mathbb{E}\int_{0}^{t}[-\big<\Lambda(\alpha_{s}) X(s),X(s)\big>
+2\big<P(\alpha_{s})b(s),X(s)\big>]ds.
\end{aligned}
$$
Set $\phi(t)\triangleq \mathbb{E}\left[|P(\alpha_{t})^{\frac{1}{2}}X(t)|^{2}\right]$, $\beta(t)\triangleq \mathbb{E}\big[|P(\alpha_{t})^{\frac{1}{2}}b(t)|^2\big]$ and $\Psi(i)\triangleq P(i)^{-\frac{1}{2}}\Lambda(i)P(i)^{-\frac{1}{2}}$.
Then
$$
\begin{aligned}
&\quad\frac{d}{dt}\phi(t)
=\mathbb{E}\big[\big<-\Psi(\alpha)P(\alpha)^{\frac{1}{2}}X+2P(\alpha)^{\frac{1}{2}}b,P(\alpha)^{\frac{1}{2}}X\big>\big]\\
&\leq -\lambda\phi+\frac{\lambda}{2}\phi+\frac{2}{\lambda}\beta
=-\frac{\lambda}{2}\phi+\frac{2}{\lambda}\beta,
\end{aligned}
$$
where $\lambda>0$ is the smallest  eigenvalue of process $\Psi(\alpha)$.
By Gronwall's inequality, we have
\begin{equation}\label{FDE-p1}
 \phi(t) \leq \phi(0)e^{-\frac{\lambda}{2}t}+\frac{2}{\lambda}\int_{0}^{t}e^{-\frac{\lambda}{2}(t-s)}\beta(s)ds.
\end{equation}
Noting that $b\in L_{\mathbb{F}}^{2}(\mathbb{R}^{n})$  implies $\int_{0}^{\infty}\beta(t)dt<\infty$.
Integrating both sides of \eqref{FDE-p1} simultaneously yields
$$
\begin{aligned}
&\quad\mathbb{E}\int_{0}^{\infty}\big<P(\alpha_{t})X(t),X(t)\big>dt=\int_{0}^{\infty}\phi(t)dt\\
&\leq\int_{0}^{\infty}\left[\phi(0)e^{-\frac{\lambda}{2}t}+\frac{2}{\lambda}\int_{0}^{t}e^{-\frac{\lambda}{2}(t-s)}\beta(s)ds\right]dt\\
&=\frac{2}{\lambda}\phi(0)+\frac{4}{\lambda^2}\int_{0}^{\infty}\beta(s)ds<\infty,
\end{aligned}
$$
which implies $X\in L_{\mathbb{F}}^{2}(\mathbb{R}^{n})$ and completes the proof.

To construct the closed-loop equilibrium strategy of Problem (M-GLQ), we need to study the solvability of the following linear BSDE:
\begin{equation}\label{BSDE}
  dY(t)=-\left[A(\alpha_{t})^{\top}Y(t)+\varphi(t)\right]dt+\mathbf{\Gamma}(t)\cdot d\mathbf{\widetilde{N}}(t), \, t\geq 0.
\end{equation}
Similar to Proposition \ref{prop-FDE}, we have the following result.

\begin{thm}\label{thm-BSDE}
Suppose Condition (A) holds and $\varphi\in L_{\mathbb{F}}^{2}(\mathbb{R}^{n})$. Then the BSDE \eqref{BSDE} admits a unique solution $\left(Y,\mathbf{\Gamma}\right)\in L_{\mathbb{F}}^{2}(\mathbb{R}^{n})\times \mathcal{D}\left(L_{\mathcal{P}}^{2}(\mathbb{R}^{n})\right)$.
\end{thm}

Before we prove it, we provide the following priori estimates, which will be used in the proof of Theorem \ref{thm-BSDE}.

\begin{prop}\label{prop-BSDE-E}
Suppose Condition (A) holds and $\varphi\in L_{\mathbb{F}}^{2}(\mathbb{R}^{n})$. Let $\left(Y,\mathbf{\Gamma}\right)\in L_{\mathbb{F}}^{2}(\mathbb{R}^{n})\times \mathcal{D}\left(L_{\mathcal{P}}^{2}(\mathbb{R}^{n})\right)$ be the solution to \eqref{BSDE}. Then
\begin{equation}\label{BSDE-E}
  \mathbb{E}\int_{0}^{\infty}\big[|Y(t)|^{2}+\sum_{j=1}^{D}\pi_{\alpha_{t}j}|\Gamma_{j}(t)|^2\big]dt
  \leq K\mathbb{E}\int_{0}^{\infty}|\varphi(t)|^{2}dt.
\end{equation}
Hereafter, $K>0$ represents a generic constant that can be different from line to line.
\end{prop}

\textbf{Proof.}
 We divide it into two steps to complete the proof. For any $t\geq 0$, we first prove that
 \begin{equation}\label{step-1}
   \mathbb{E}\!\Big\{\!|Y(t)|^{2}\!+\!\sum_{j=1}^{D}\int_{t}^{\infty}\pi_{\alpha_{s}j}|\Gamma_{j}(s)|^2ds\!\Big\}\!\leq\! K\mathbb{E}\!\int_{t}^{\infty}\!|\varphi(s)|^{2}\!ds,
 \end{equation}
 and then, we verify that
  \begin{equation}\label{step-2}
   \mathbb{E}\int_{0}^{\infty}\!|Y(t)|^{2}\!dt\leq K\mathbb{E}\int_{0}^{\infty}\!\big[|\varphi(t)|^{2}\!+\!\sum_{j=1}^{D}\pi_{\alpha_{t}j}|\Gamma_{j}(t)|^2\big]\!dt.
 \end{equation}
 The desired result \eqref{BSDE-E} follows from \eqref{step-1} and \eqref{step-2} directly. Now, let us concretize the above two steps. Since Condition (A) holds, let  $\mathbf{P}\in \mathfrak{P}$. Then we can choose $\varepsilon>0$ such that
 $$\Lambda^{\varepsilon}(i)\!\triangleq \!-\big[\!P(i)A(i)+A(i)^{\top}P(i)\!+\!(1\!+\!\varepsilon)\sum_{j=1}^{D}\pi_{ij}P(j)\!\big]>0,$$
Applying It\^o's rule to $\big<P\left(\alpha_{t}\right)^{-1}Y(t),Y(t)\big>$, we obtain
{ \scriptsize
 \begin{align*}
 &\quad\mathbb{E}\Big[\big<P\left(\alpha_{T}\right)^{-1}Y(T),Y(T)\big>-\big<P\left(\alpha_{t}\right)^{-1}Y(t),Y(t)\big>\Big]\\
 &=\mathbb{E}\int_{t}^{T}\Big\{-\big<\big[A(\alpha)^{\top}P(\alpha)+P(\alpha)A(\alpha)\big]P(\alpha)^{-1}Y,P(\alpha)^{-1}Y\big>\\
 &\quad+\sum_{j=1}^{D}\pi_{\alpha_{s}j}\big<P(j)^{-1}\Gamma_{j},\Gamma_{j}\big>-2\sum_{j=1}^{D}\pi_{\alpha_{s}j}\big<\Gamma_{j},P(\alpha)^{-1}Y\big>\\
 &\quad+\sum_{j=1}^{D}\pi_{\alpha_{s}j}\big<P(j)^{-1}Y,Y\big>+2\sum_{j=1}^{D}\pi_{\alpha_{s}j}\big<P(j)^{-1}\Gamma_{j},Y\big>\\
 &\quad-2\big<\varphi,P\left(\alpha\right)^{-1}Y\big>\Big\}ds\\
 &\geq\mathbb{E}\int_{t}^{T}\Big\{\big<\Lambda^{\varepsilon}(\alpha)P(\alpha)^{-1}Y,P(\alpha)^{-1}Y\big>-2\big<\varphi,P(\alpha)^{-1}Y\big>\\
 &\quad+\sum_{j=1}^{D}\pi_{\alpha_{s}j}\left[(1+\varepsilon)\big<P(j)P(\alpha)^{-1}Y,P(\alpha)^{-1}Y\big>-2\big<\Gamma_{j},P(\alpha)^{-1}Y\big>\right]\\
 &\quad+\sum_{j=1}^{D}\pi_{\alpha_{s}j}\big<P(j)^{-1}\Gamma_{j},\Gamma_{j}\big>+2\sum_{j=1}^{D}\pi_{\alpha_{s}j}\big<P(j)^{-1}\Gamma_{j},Y\big>\Big\}ds\\
 &=\mathbb{E}\int_{t}^{T}\Big\{\big<\Lambda^{\varepsilon}(\alpha)P(\alpha)^{-1}Y,P(\alpha)^{-1}Y\big>-2\big<\varphi,P(\alpha)^{-1}Y\big>\\
 &\quad+(1+\varepsilon)\sum_{j=1}^{D}\pi_{\alpha_{s}j}\left|P(j)^{\frac{1}{2}}\left[\!P(\alpha)^{-1}Y\!-\!\frac{1}{1+\varepsilon}P(j)^{-1}\Gamma_{j}\!\right]\right|^{2}\\
 &\quad+\frac{\varepsilon}{1+\varepsilon}\sum_{j=1}^{D}\pi_{\alpha_{s}j}\big<P(j)^{-1}\Gamma_{j},\Gamma_{j}\big>+2\sum_{j=1}^{D}\pi_{\alpha_{s}j}\big<P(j)^{-1}\Gamma_{j},Y\big>\Big\}ds
 \end{align*}
}%
Let $\lambda>0$ be the smallest  eigenvalue of process $\Lambda(\alpha)$ and $\mu>0 \text{ }(\sqrt{\beta}>0)$ be the smallest (largest) eigenvalue of $P(\alpha)^{-1}$.
 Then
 \begin{align*}
 &\quad\mathbb{E}\left[\left<P\left(\alpha_{T}\right)^{-1}Y(T),Y(T)\right>-\left<P\left(\alpha_{t}\right)^{-1}Y(t),Y(t)\right>\right]\\
 &\geq \mathbb{E}\int_{t}^{T}\Big\{\lambda\left|P(\alpha)^{-1}Y\right|^{2}-\lambda\left|P(\alpha)^{-1}Y\right|^{2}-\frac{1}{\lambda}\left|\varphi\right|^{2}\\
 &\quad+\frac{\varepsilon\mu}{1+\varepsilon}\sum_{j=1}^{D}\pi_{\alpha_{s}j}\left|\Gamma_{j}\right|^{2}-\sum_{j=1}^{D}\pi_{\alpha_{s}j}\left[\theta\beta\left|\Gamma_{j}\right|^{2}+\frac{1}{\theta}\left|Y\right|^{2}\right]\Big\}ds\\
 &=\mathbb{E}\int_{t}^{T}\Big\{-\frac{1}{\lambda}\left|\varphi\right|^{2}+\left[\frac{\varepsilon\mu}{1+\varepsilon}-\theta\beta\right]\sum_{j=1}^{D}\pi_{\alpha_{s}j}\left|\Gamma_{j}\right|^{2}\Big\}ds.
 \end{align*}
 Set $\theta=\frac{\varepsilon\mu}{2\beta(1+\varepsilon)}>0$. The above equation implies
 {\scriptsize
 \begin{align*}
&\quad\mathbb{E}\Big\{\left<P\left(\alpha_{t}\right)^{-1}Y(t),Y(t)\right>
+\frac{\varepsilon\mu}{2(1+\varepsilon)}\int_{t}^{T}\sum_{j=1}^{D}\pi_{\alpha_{s}j}\left|\Gamma_{j}\right|^{2}ds\Big\}\\
 &\leq \mathbb{E}\Big\{\left<P\left(\alpha_{T}\right)^{-1}Y(T),Y(T)\right>
 +\int_{t}^{T}\frac{1}{\lambda}\left|\varphi(s)\right|^{2}ds\Big\}.
 \end{align*}
}%
 Since $Y\in L_{\mathbb{F}}^{2}(\mathbb{R}^{n})$, we must have $\lim_{T\rightarrow\infty}\mathbb{E}\left|Y(T)\right|^{2}=0$. 
 The result \eqref{step-1} follows by letting $T\rightarrow\infty$.

 For step $2$, applying It\^o's rule to $\big<P\left(\alpha_{t}\right)^{-1}Y(t),Y(t)\big>$ again, we obtain
 {\scriptsize
 \begin{align*}
&\quad\mathbb{E}\left[\left<P\left(\alpha_{T}\right)^{-1}Y(T),Y(T)\right>-\left<P\left(i\right)^{-1}Y(0),Y(0)\right>\right]\\
&=\mathbb{E}\int_{0}^{T}\Big\{-\big<\big[A(\alpha)^{\top}P(\alpha)+P(\alpha)A(\alpha)\big]P(\alpha)^{-1}Y,P(\alpha)^{-1}Y\big>\\
&\quad+\sum_{j=1}^{D}\pi_{\alpha_{s}j}\big[\big<P(j)^{-1}\left[Y+\Gamma_{j}\right],Y+\Gamma_{j}\big>-2\big<\Gamma_{j},P(\alpha)^{-1}Y\big>\big]\\
 &\quad-2\big<\varphi,P(\alpha)^{-1}Y\big>\Big\}ds.
 \end{align*}
}
 Let $\mu>0$ be the smallest eigenvalue of the process $-\left[A(\alpha)^{\top}P(\alpha)+P(\alpha)A(\alpha)\right]$. Then
 \begin{align*}
&\quad\mathbb{E}\left[\left<P\left(\alpha_{T}\right)^{-1}Y(T),Y(T)\right>-\left<P\left(i\right)^{-1}Y(0),Y(0)\right>\right]\\
&\geq \mathbb{E}\int_{0}^{T}\Big\{\mu\left|P(\alpha)^{-1}Y\right|^{2}-\frac{\mu}{2}\left|P(\alpha)^{-1}Y\right|^{2}-\frac{2}{\mu}\left|\varphi\right|^{2}\\
&\quad-\sum_{j=1}^{D}\pi_{\alpha_{s}j}\left|\Gamma_{j}\right|^{2}
 \Big\}ds\\
 &=\mathbb{E}\int_{0}^{T}\Big\{\frac{\mu}{2}\left|P(\alpha)^{-1}Y\right|^{2}-\frac{2}{\mu}\left|\varphi\right|^{2}-\sum_{j=1}^{D}\pi_{\alpha_{s}j}\left|\Gamma_{j}\right|^{2}\Big\}ds.
 \end{align*}
Consequently,
\begin{align*}
&\mathbb{E}\int_{0}^{T}\frac{\mu}{2}\left|P(\alpha)^{-1}Y\right|^{2}\leq\mathbb{E}\Big\{\big<P(\alpha_{T})^{-1}Y(T),Y(T)\big>\\
&\qquad+\int_{0}^{T}\Big[
\frac{2}{\mu}\left|\varphi\right|^{2}+\sum_{j=1}^{D}\pi_{\alpha_{s}j}\left|\Gamma_{j}\right|^{2}\Big]ds
\Big\}.
\end{align*}
Letting $T\rightarrow\infty$, the \eqref{step-2} can be derived immediately. This completes the proof.

\textbf{Proof of Theorem \ref{thm-BSDE}.}
  The uniqueness follows immediately from priori estimate \eqref{BSDE-E}. For the existence,  let
  \begin{align*}
    &\varphi^{(k)}(t)=\varphi(t)\mathbb{I}_{[0,k]}(t), \qquad\qquad t\in[0,\infty), \\
   & \widehat{\varphi}^{(k)}(t)=\varphi(t),  \qquad\qquad\qquad\quad  t\in[0,k].
  \end{align*}
  Then $\left\{\varphi^{(k)}\right\}_{k=1}^{\infty}$ converges to $\varphi$ in $L_{\mathbb{F}}^{2}(\mathbb{R}^{n})$. Now, consider the following finite horizon BSDE
  {\scriptsize
  \begin{equation}\label{BSDE-F}
    \left\{
    \begin{aligned}
    &d\widehat{Y}^{(k)}(t)=-\left[A(\alpha_{t})^{\top}\widehat{Y}^{(k)}(t)+\widehat{\varphi}^{(k)}(t)\right]dt
    +\mathbf{\widehat{\Gamma}^{(k)}}(t)\cdot d\mathbf{\widetilde{N}}(t), \\
    &\widehat{Y}^{(k)}(k)=0,\quad t\in[0,k].
    \end{aligned}
    \right.
  \end{equation}
  }%
  Obviously, for any $k>0$, the above BSDE admits a square integrable solution $\big(\widehat{Y}^{(k)},\mathbf{\widehat{\Gamma}^{(k)}}\big)$.
  Setting
  \begin{equation*}\left\{
  \begin{aligned}
  &Y^{(k)}(t)=\widehat{Y}^{(k)}(t)\mathbb{I}_{[0,k]}(t)+0\mathbb{I}_{[k,\infty)}(t),\\
  &\mathbf{\Gamma^{(k)}}(t)=\mathbf{\widehat{\Gamma}^{(k)}}(t)\mathbb{I}_{[0,k]}(t)+\mathbf{0}\mathbb{I}_{[k,\infty)}(t), \quad t\in [0,\infty),
  \end{aligned}
  \right.
  \end{equation*}
  then $\left(Y^{(k)},\mathbf{\Gamma^{(k)}}\right)\in L_{\mathbb{F}}^{2}(\mathbb{R}^{n})\times \mathcal{D}\left(L_{\mathcal{P}}^{2}(\mathbb{R}^{n})\right)$ solves  BSDE:
  $$
  dY^{(k)}(t)\!=\!-\!\left[A(\alpha_{t})^{\top}Y^{(k)}(t)+\varphi^{(k)}(t)\right]dt+\mathbf{\Gamma^{(k)}}(t)\cdot d\mathbf{\widetilde{N}}(t).
  $$
On the other hand, by Proposition \ref{prop-BSDE-E}, for any $k$, $l>0$, we have
    \begin{align*}
    &\mathbb{E}\int_{0}^{\infty}\Big[\left|Y^{(k)}-Y^{(l)}\right|^{2}
    +\sum_{j=1}^{D}\pi_{\alpha_{t}j}\left|\Gamma_{j}^{(k)}-\Gamma_{j}^{(l)}\right|^{2}\Big]dt\\
    \leq &K\mathbb{E}\int_{0}^{\infty}\left|\varphi_{j}^{(k)}-\varphi_{j}^{(l)}\right|^{2}dt.
    \end{align*}
    Therefore, there exists a triple $\left(Y,\mathbf{\Gamma}\right)\in L_{\mathbb{F}}^{2}(\mathbb{R}^{n})\times\mathcal{D}\left(L_{\mathcal{P}}^{2}(\mathbb{R}^{n})\right)$ such that
    \begin{align*}
    \quad\mathbb{E}\int_{0}^{\infty}\Big[\left|Y^{(k)}-Y\right|^{2}
    +\sum_{j=1}^{D}\pi_{\alpha_{t}j}\left|\Gamma_{j}^{(k)}-\Gamma\right|^{2}\Big]dt\rightarrow 0,
    \end{align*}
   as $k\rightarrow\infty$. One can easily verify that $\left(Y,\mathbf{\Gamma}\right)$ is the unique solution of BSDE \eqref{BSDE}. This completes the proof.


\section{The LQ problem on infinite horizon}\label{section-LQ}
This section studies the  LQ problem with Markovian jumps in infinite horizon. The results obtained will be used to solve the Problem (M-GLQ).
Consider the following controlled state equation:
\begin{equation}\label{state-LQ}
  \left\{
 \begin{aligned}
   \dot{X}(t)&=A\left(\alpha_{t}\right)X(t)+B\left(\alpha_{t}\right)u(t)+b(t),\\
   X(0)&=x,\quad \alpha(0)=i,
   \end{aligned}
  \right.
\end{equation}
with cost functional:
\begin{equation}\label{cost-LQ}
\begin{aligned}
    J\left(x,i;u\right)
    & \triangleq \mathbb{E}\int_{0}^{\infty}\left[\!
    \left<\!
    \left(\!
    \begin{matrix}
    Q(\alpha)&S(\alpha)^{\top} \\
    S(\alpha)&R(\alpha)
    \end{matrix}
    \!\right)\!
    \left(\!
    \begin{matrix}
    X\\
    u
    \end{matrix}
    \!\right)\!,
    \left(\!
    \begin{matrix}
    X\\
    u
    \end{matrix}
    \!\right)
    \!\right>\!\right.\\
    &\quad\left.+2
    \left<\!
    \left(\!
    \begin{matrix}
    q\\
    \rho
    \end{matrix}
    \!\right)\!,
    \left(\!
    \begin{matrix}
    X\\
    u
    \end{matrix}
    \!\right)
    \!\right>\!\right]dt.
  \end{aligned}
\end{equation}
Here, $b,\, q\in L_{\mathbb{F}}^{2}(\mathbb{R}^{n})$ and $\rho\in L_{\mathbb{F}}^{2}(\mathbb{R}^{m})$. For fixed $i\in\mathcal{S}$, $A(i)\in\mathbb{R}^{n\times n}$, $B(i),\text{} S(i)^{\top}\in\mathbb{R}^{n\times m}$, $Q(i)\in\mathbb{S}^{n}$, $R(i)\in\mathbb{S}_{+}^{m}$ and $Q(i)-S(i)^{\top}R(i)^{-1}S(i)\geq 0$. The solution of \eqref{state-LQ} is denoted by $X\equiv X(\cdot;x,i,u)$.
For any given $(x,i)\in \mathbb{R}^{n}\times \mathcal{S}$, we define the admissible control sets of LQ problem  as follows:
$$\mathcal{U}_{LQ}(x,i)\triangleq\left\{u\in L_{\mathcal{P}}^{2}(\mathbb{R}^{m}) \mid X(\cdot;x,i,u)\in L_{\mathbb{F}}^{2}(\mathbb{R}^{n})\right\}. $$
\textbf{Problem (M-LQ).} For any given $(x,i)\in \mathbb{R}^{n}\times \mathcal{S}$, find a $u^{*}\in\mathcal{U}_{LQ}(x,i)$ such that
\begin{equation}\label{value-LQ}
    J\left(x,i;u^{*}\right)=\inf_{u\in\mathcal{U}_{LQ}(x,i)}J\left(x,i;u\right)\triangleq V(x,i) .
\end{equation}
Any $u^{*}\in\mathcal{U}_{LQ}(x,i)$ satisfying \eqref{value-LQ} is called an open-loop optimal control of Problem (M-LQ) for the initial value $(x,i)$, and the corresponding $X^{*}\equiv X(\cdot;x,i,u^{*})$ is called an optimal state process. The function $V(\cdot,\cdot)$ is called the value function of Problem (M-LQ). In addition, if $b,\,q,\,\rho=0$, then corresponding state process, cost functional, value function and problem are denoted by $X^{0}\equiv X^{0}(\cdot;x,i,u)$, $J^{0}(x,i;u)$, $V^{0}(x,i)$ and Problem (M-LQ)$^{0}$, respectively.

\begin{defn}\label{def-stabilizable}
The system \eqref{state-LQ} is said to be $L^{2}$-stabilizable if and only if there exists a $\mathbf{\Theta}\in\mathcal{D}\big(\mathbb{R}^{m\times n}\big)$  such that for any $i\in\mathcal{S}$,
{\scriptsize
\begin{equation}\label{L-2}
P(i)[A(i)+B(i)\Theta(i)]+[A(i)+B(i)\Theta(i)]^{\top}P(i)+\sum_{j=1}^{D}\pi_{ij}P(j)<0,
\end{equation}
}
for some $\mathbf{P}\in \mathcal{D}\left(\mathbb{S}_{+}^{n}\right)$. Any element $\mathbf{\Theta}\in\mathcal{D}\big(\mathbb{R}^{m\times n}\big)$ satisfying the above condition is called a stabilizer of system \eqref{state-LQ} and we denote $\mathcal{H}_{\alpha}$ represents the set consisting of all stabilizer.
\end{defn}
We will be working with the following assumption in the rest of this section.

\textbf{(H1)} System  \eqref{state-LQ} is $L^{2}$-stabilizable, i.e, $\mathcal{H}_{\alpha}\neq\emptyset$.

To simplify our further analysis, for any given $\mathbf{\Theta}\in \mathcal{H}_{\alpha}$, we consider the following state equation:
\begin{equation}\label{state-LQ-Theta-nu}
       \left\{
       \begin{aligned}
     dX_{\Theta}(t)&=\left\{A_{\Theta}(\alpha_{t})X_{\Theta}(t)+B\left(\alpha_{t}\right)\nu(t)+b(t)\right\}dt,\\
     X_{\Theta}(0)&=x,\quad \alpha(0)=i,\quad t\geq0,
       \end{aligned}
       \right.
  \end{equation}
  and cost functional
  \begin{equation}\label{cost-LQ-Theta-nu}
\begin{aligned}
&J_{\Theta}(x,i;\nu)\triangleq J(x,i;\Theta(\alpha)X_{\Theta}+\nu)\\
&= \mathbb{E}\int_{0}^{\infty}\left[\!
    \left<\!
    \left(\!
  \begin{matrix}
   Q_{\Theta}(\alpha)&S_{\Theta}(\alpha)^{\top}\\
   S_{\Theta}(\alpha)&R(\alpha)
 \end{matrix}
    \!\right)\!
    \left(\!
\begin{matrix}
    X_{\Theta}\\
    \nu
\end{matrix}
   \! \right)\!,
    \left(\!
  \begin{matrix}
     X_{\Theta}\\
    \nu
\end{matrix}
    \!\right)\!
    \right>\right.\\
  &\quad \left.+2\left<\!
    \left(\!
    \begin{matrix}
 q_{\Theta}\\
    \rho
\end{matrix}
    \!\right)\!,
    \left(\!
   \begin{matrix}
     X_{\Theta}\\
    \nu
\end{matrix}
   \! \right)\!
    \right>\!\right]dt,
  \end{aligned}
\end{equation}
where
\begin{equation*}
\begin{aligned}
A_{\Theta}(i)&=A(i)+B(i)\Theta(i),\quad S_{\Theta}(i)=S(i)+R(i)\Theta(i),\\
Q_{\Theta}(i)&=Q(i)\!+\!S(i)^{\top}\Theta(i)\!+\!\Theta(i)^{\top}S(i)\!+\!\Theta(i)^{\top}R(i)\Theta(i), \\
q_{\Theta}&=q+\Theta(\alpha\cdot)^{\top}\rho.
\end{aligned}
\end{equation*}
Obviously,  by Proposition \ref{prop-FDE}, the equation \eqref{state-LQ-Theta-nu} admits a unique solution  $X_{\Theta}\equiv X_{\Theta}(\cdot;x,i,\nu)\in L_{\mathbb{F}}^{2}(\mathbb{R}^{n})$ for any $\nu\in L_{\mathcal{P}}^{2}(\mathbb{R}^{m})$.

\textbf{Problem (M-LQ)$_{\Theta}$.} For any given $(x,i)\in \mathbb{R}^{n}\times \mathcal{S}$, find a $\nu^{*}\in L_{\mathcal{P}}^{2}(\mathbb{R}^{m})$ such that
\begin{equation}\label{value-LQ-Theta}
    J_{\Theta}\left(x,i;\nu^{*}\right)=\inf_{\nu\in L_{\mathcal{P}}^{2}(\mathbb{R}^{m})}J_{\Theta}\left(x,i;\nu\right).
\end{equation}
Any $\nu^{*}\in L_{\mathcal{P}}^{2}(\mathbb{R}^{m})$ satisfying \eqref{value-LQ-Theta} is called an open-loop optimal control of Problem (M-LQ)$_{\Theta}$ for the initial value $(x,i)$. In addition, if $b$, $q_{\Theta}$, $\rho=0$, we further denote the corresponding state process, cost functional and problem as $X_{\Theta}^{0}\equiv X_{\Theta}^{0}(\cdot;x,i,\nu)$, $J_{\Theta}^{0}(x,i;\nu)$ and Problem (M-LQ)$_{\Theta}^{0}$, respectively.

\begin{defn}
\begin{description}
  \item[(i)] Problem (M-LQ) is said to be well-posedness if
  $$
  V(x,i)>-\infty,\qquad \forall (x,i)\in\mathbb{R}^{n}\times\mathcal{S}.
  $$
  \item[(ii)] A pair $(\widehat{\mathbf{\Theta}},\widehat{\nu})\in\mathcal{H}_{\alpha}\times L_{\mathcal{P}}^{2}(\mathbb{R}^{m})$ is called a closed-loop optimal control of Problem (M-LQ) if for any $\left(x,i\right)\in \mathbb{R}^{n}\times\mathcal{S}$ and $ u\in \mathcal{U}_{LQ}(x,i)$, the following holds:
  \begin{equation}\label{closed}
   J(x,i;\widehat{\Theta}(\alpha)X_{\widehat{\Theta}}+\widehat{\nu})\leq J(x,i;u).
  \end{equation}
  In addition, Problem (M-LQ) is said to be (uniquely) closed-loop solvable if it admits a (unique) closed-loop optimal control.
\end{description}
\end{defn}

The following proposition provides a profound characterization of the closed-loop solvability for Problem (M-LQ).

\begin{prop}\label{prop-LQ-control-characterization}
Let assumption (H1) hold. A pair $(\widehat{\mathbf{\Theta}},\widehat{\nu})\in\mathcal{H}_{\alpha}\times L_{\mathcal{P}}^{2}(\mathbb{R}^{m})$ is a closed-loop optimal control of Problem (M-LQ) if and only if $\widehat{\nu}$ is an open-loop optimal control of Problem (M-LQ)$_{\widehat{\Theta}}$ for any initial value $(x,i)\in\mathbb{R}^{n}\times\mathcal{S}$.
\end{prop}
\textbf{Proof.}
By definition, we only need to prove that
 \begin{align*}
\mathcal{U}_{\widehat{\Theta}}(x,i)&\triangleq\left\{\widehat{\Theta}(\alpha)X_{\widehat{\Theta}}(\cdot;x,i,\nu)+\nu\Big|\nu\in L_{\mathcal{P}}^{2}(\mathbb{R}^{m})\right\}\\
&=\mathcal{U}_{LQ}(x,i),
 \end{align*}
 holds for any $(x,i)\in\mathbb{R}^{n}\times\mathcal{S}$.
It follows from  Proposition \ref{prop-FDE} that $X_{\widehat{\Theta}}(\cdot;x,i,\nu)\in L_{\mathbb{F}}^{2}(\mathbb{R}^{n})$ for any $\nu\in L_{\mathcal{P}}^{2}(\mathbb{R}^{m})$. Consequently, we have
  $$
  u\triangleq\widehat{\Theta}(\alpha)X_{\widehat{\Theta}}(\cdot;x,i,\nu)+\nu\in L_{\mathcal{P}}^{2}(\mathbb{R}^{m}).
  $$
  By uniqueness, one has $X_{\widehat{\Theta}}(\cdot;x,i,\nu)=X(\cdot;x,i,u)$, which implies $u\in \mathcal{U}_{LQ}(x,i)$ and $\mathcal{U}_{\widehat{\Theta}}(x,i)\subseteq\mathcal{U}_{LQ}(x,i)$.
  On the other hand, for any $u\in\mathcal{U}_{LQ}(x,i)$, one has
  $$
   \nu\triangleq u-\widehat{\Theta}(\alpha)X(\cdot;x,i,u)\in L_{\mathcal{P}}^{2}(\mathbb{R}^{m}).
     $$
  By uniqueness again, we have $X(\cdot;x,i,u)=X_{\widehat{\Theta}}(\cdot;x,i,\nu)$. Consequently, $u\in \mathcal{U}_{\widehat{\Theta}}(x,i)$ and $\mathcal{U}_{LQ}(x,i)\subseteq\mathcal{U}_{\widehat{\Theta}}(x,i)$.
  In summary, the proposition is proven.

Next, we consider the adjoint equation of \eqref{state-LQ-Theta-nu} in Problem (M-LQ)$_{\Theta}$ as follows:
\begin{equation}\label{LQ-BSDE-Theta-nu}
\begin{aligned}
 dY_{\Theta}(t)&=-\big\{A_{\Theta}(\alpha)^{\top}Y_{\Theta}+Q_{\Theta}(\alpha)X_{\Theta}+S_{\Theta}(\alpha)^{\top}\nu\\
 &\quad+q_{\Theta}\big\}dt+\mathbf{\Gamma}_{\Theta}(t)\cdot d\mathbf{\widetilde{N}}(t),\quad t\geq 0.
\end{aligned}
\end{equation}
Obviously, for any $(\mathbf{\Theta},\nu)\in\mathcal{H}_{\alpha}\times L_{\mathcal{P}}^{2}(\mathbb{R}^{m})$, the BSDE \eqref{LQ-BSDE-Theta-nu} admits a unique $L^2$-integrable solution  and we denote it as
$$
\big(Y_{\Theta},\mathbf{\Gamma}_{\Theta}\big)\equiv\big(Y_{\Theta}(\cdot;x,i,\nu),\mathbf{\Gamma}_{\Theta}(\cdot;x,i,\nu)\big).
$$
Additionally, if $b$, $q$ $\rho=0$,
then we denote the corresponding solutions to \eqref{state-LQ-Theta-nu} and \eqref{LQ-BSDE-Theta-nu}  as
$$ \big(X_{\Theta}^{0},Y_{\Theta}^{0},\mathbf{\Gamma}_{\Theta}^{0}\big)\equiv\big(X_{\Theta}^{0}(\cdot;x,i,\nu),Y_{\Theta}^{0}(\cdot;x,i,\nu),\mathbf{\Gamma}_{\Theta}^{0}(\cdot;x,i,\nu)\big).$$
With the above notations, we have the following result.

\begin{thm}\label{thm-SLQ-FBSDE}
Suppose assumption (H1) holds. A pair $(\widehat{\mathbf{\Theta}},\widehat{\nu})\in\mathcal{H}_{\alpha}\times L_{\mathcal{P}}^{2}(\mathbb{R}^{m})$ is a closed-loop optimal control of Problem (M-LQ) if and only if, for any $(x,i)\in\mathbb{R}^{n}\times\mathcal{S}$, the following stationarity condition holds:
\begin{equation}\label{SLQ-stationarity}
   \begin{aligned}
      B(\alpha)^{\top}Y_{\widehat{\Theta}}\!+\!S_{\widehat{\Theta}}(\alpha)X_{\widehat{\Theta}}\!+\!R(\alpha)\widehat{\nu}\!+\!\rho=0.
   \end{aligned}
\end{equation}
\end{thm}

\textbf{Proof.}
The proof is conventional (see, for example, \cite{Sun.J.R.2016_open-closed,zhang2021open}). We leave the detailed derivations to the interested readers.

 \begin{rem}\label{rmk-SLQ-Y}
Suppose $(\widehat{\mathbf{\Theta}},\widehat{\nu})\in\mathcal{H}_{\alpha}\times L_{\mathcal{P}}^{2}(\mathbb{R}^{m})$ is a closed-loop optimal control of Problem (M-LQ). Using the stationarity condition \eqref{SLQ-stationarity}, one can further simplify the  adjoint equation \eqref{LQ-BSDE-Theta-nu} as follows:
\begin{align*}
d\widehat{Y}(t)&=-\big\{A(\alpha)^{\top}\widehat{Y}+\big(Q(\alpha)+S(\alpha)^{\top}\widehat{\Theta}(\alpha)\big)\widehat{X}\\
&\quad+S(\alpha)^{\top}\widehat{\nu}(t)+q(t)\big\}dt+\widehat{\mathbf{\Gamma}}(t)\cdot d\mathbf{\widetilde{N}}(t),
\end{align*}
where
$$
\left(\widehat{X},\widehat{Y},\widehat{\mathbf{\Gamma}}\right)\equiv
\left(X_{\widehat{\Theta}}(\cdot;x,i,\widehat{\nu}),Y_{\widehat{\Theta}}(\cdot;x,i,\widehat{\nu}),\mathbf{\Gamma}_{\widehat{\Theta}}(\cdot;x,i,\widehat{\nu})\right).
$$
 \end{rem}

We have the following corollary to Theorem \ref{thm-SLQ-FBSDE}.

\begin{cor}\label{coro-SLQ-0}
If $(\widehat{\mathbf{\Theta}},\widehat{\nu})\in\mathcal{H}_{\alpha}\times L_{\mathcal{P}}^{2}(\mathbb{R}^{m})$  is a closed-loop optimal control of Problem (M-LQ), then $(\widehat{\mathbf{\Theta}},0)$ is a closed-loop optimal control of Problem (M-LQ)$^{0}$.
\end{cor}

\textbf{Proof.}
By Theorem \ref{thm-SLQ-FBSDE}, we just need to verify that, for any $(x,i)\in\mathbb{R}^{n}\times\mathcal{S}$, the following  holds:
\begin{equation}\label{SLQ-stationarity-0}
   \begin{aligned}
      B(\alpha)^{\top}Y_{\widehat{\Theta}}^{0}(\cdot;x,i,\widehat{\nu})
      +S_{\widehat{\Theta}}(\alpha)X_{\widehat{\Theta}}^{0}(\cdot;x,i,\widehat{\nu})=0.
   \end{aligned}
\end{equation}
Since the closed-loop control $(\widehat{\mathbf{\Theta}},\widehat{\nu})$ and the condition \eqref{SLQ-stationarity} are independent of initial value $(x,i)\in \mathbb{R}^{n}\times\mathcal{S}$, we have
\begin{equation}\label{SLQ-0-P1}
\begin{aligned}
     &B(\alpha)^{\top}Y_{\widehat{\Theta}}(\cdot;0,i,\widehat{\nu})+S_{\widehat{\Theta}}(\alpha)X_{\widehat{\Theta}}(\cdot;0,i,\widehat{\nu})\\
      &\quad+R(\alpha_{t})\widehat{\nu}+\rho=0.
\end{aligned}
\end{equation}
The desired result \eqref{SLQ-stationarity-0} follows by subtracting \eqref{SLQ-0-P1} from \eqref{SLQ-stationarity}.

To simplify the notations, for any given $\mathbf{P}\in\mathcal{D}\left(\mathbb{S}^{n}\right)$, let
$$
    \begin{array}{l}
    \mathcal{M}(P,i)\triangleq P(i)A(i)\!+\!A(i)^{\top}P(i)\!+\!Q(i)\!+\!\sum_{j=1}^{D}\!\pi_{ij}P(j),\\
    \mathcal{L}(P,i)\triangleq P(i)B(i)+S(i)^{\top},\\
    \widehat{q}=P(\alpha)b(t)+q.
    \end{array}
$$
The following lemma provides another representation of cost functional for Problem (M-LQ).
 \begin{lem}\label{lem-useful}
For any given $\mathbf{P}\in\mathcal{D}\left(\mathbb{S}^{n}\right)$, the cost functional \eqref{cost-LQ} admits the following representation:
\begin{align*}
  &J(x,i;u)\!=\!\mathbb{E}\!\int_{0}^{\infty}\!\left[\!
    \left<\!
    \left(\!
    \begin{matrix}
    \mathcal{M}(P,\alpha)&\mathcal{L}(P,\alpha) \\
    \mathcal{L}(P,\alpha)^{\top}&R(\alpha)
    \end{matrix}
    \!\right)\!
    \left(\!
    \begin{matrix}
    X\\
    u
    \end{matrix}
    \!\right)\!,
    \left(\!
    \begin{matrix}
    X\\
    u
    \end{matrix}
    \!\right)\!
    \right>\!\right.\\
   &\quad\left. +2\left<\!\left(\!\begin{matrix}
    \widehat{q}\\
    \rho
    \end{matrix}
    \!\right)\!,
    \left(\!
    \begin{matrix}
    X\\
    u
    \end{matrix}
    \!\right)\!
    \right>
    \!\right]dt+\big<P(i)x,x\big>.
    \end{align*}
 \end{lem}
\textbf{Proof.} This is immediate by applying It\^o's rule to $\left<P(\alpha_{t})X(t),X(t)\right>$.

\begin{rem}\label{rem-useful}
If $b,\,q,\,\rho=0$, then  obviously we have
 \begin{align*}
    &J^{0}(x,i;u)\!=\!\big<P(i)x,x\big>\\
  &+\mathbb{E}\!\int_{0}^{\infty}\!\left[\!
    \left<\!
    \left(\!
    \begin{matrix}
    \mathcal{M}(P,\alpha)&\mathcal{L}(P,\alpha) \\
    \mathcal{L}(P,\alpha)^{\top}&R(\alpha)
    \end{matrix}
    \!\right)\!
    \left(\!
    \begin{matrix}
    X^{0}\\
    u
    \end{matrix}
    \!\right)\!,
    \left(\!
    \begin{matrix}
    X^{0}\\
    u
    \end{matrix}
    \!\right)\!
    \right>\!\right]dt.
  \end{align*}
\end{rem}

The following technical result can be shown by a simple adaptation of the well-known result in a deterministic case (see, e.g., Anderson and Moore \cite{Anderson.1990})
 \begin{lem}\label{lem-finite-1}
Problem (M-LQ)$^{0}$ is well-posedness if and only if there exists a $\mathbf{P}\in\mathcal{D}\left(\mathbb{S}^{n}\right)$ such that
$$V^{0}(x,i)=\big<P(i)x,x\big>,\qquad \forall (x,i)\in\mathbb{R}^{n}\times\mathcal{S}.$$
 \end{lem}
Now, we define the following convex set
\begin{equation}\label{hua-G}
  \mathcal{G}\triangleq\left\{\mathbf{P}\in\mathcal{D}\left(\mathbb{S}^{n}\right)\Big|\left[
\begin{matrix}
  \mathcal{M}(P,i)  &  \mathcal{L}(P,i)  \\
  \mathcal{L}(P,i)^{\top}  & R(i)
\end{matrix}
\right]\geq 0,\,\forall i\in\mathcal{S}\right\}.
\end{equation}
An element $\mathbf{P}\in\mathcal{G}$ is called the maximal element if and only if, for any $\mathbf{\widetilde{P}}\in\mathcal{G}$, we have $P(i)\geq\widetilde{P}(i)$ for all $i\in \mathcal{S}$.

\begin{prop}\label{prop-finite-2}
 Problem (M-LQ)$^{0}$ is well-posedness if and only if $\mathcal{G}\neq\emptyset$. In this case, $\mathcal{G}$ has a unique maximal element $\mathbf{P}\in\mathcal{G}$ such that
 $$
 V^{0}(x,i)=\big<P(i)x,x\big>,\qquad \forall (x,i)\in\mathbb{R}^{n}\times\mathcal{S}.
 $$
\end{prop}
 \textbf{Proof.}
 Let $\mathbf{\widetilde{P}}\in\mathcal{G}$. By Remark \ref{rem-useful}, we have
   \begin{equation}\label{proof-finite-2-1}
\begin{aligned}
   &\quad J^{0}(x,i;u)\\
   &=\mathbb{E}
    \int_{0}^{\infty}\!
    \left<\!
    \left(\!
    \begin{matrix}
    \mathcal{M}(\widetilde{P},\alpha)&\mathcal{L}(\widetilde{P},\alpha)\\
    \mathcal{L}(\widetilde{P},\alpha)^{\top}&R(\alpha)
    \end{matrix}
    \!\right)\!
    \left(\!
    \begin{matrix}
    X^{0}\\
    u
    \end{matrix}
    \!\right)\!,
    \!\left(\!
    \begin{matrix}
    X^{0}\\
    u
    \end{matrix}
    \!\right)\!\right>dt\\
    &\quad+\big<\widetilde{P}(i)x,x\big>\\
    &\geq \big<\widetilde{P}(i)x,x\big>.
  \end{aligned}
\end{equation}
Consequently,
$$V^{0}(x,i)\!=\!\inf_{u\in\mathcal{U}_{ad}(x,i)}J^{0}(x,i;u)\geq \big<\widetilde{P}(i)x,x\big>>-\infty.$$
Conversely, if problem (M-LQ)$^{0}$ is well-posedness, then by Lemma \ref{lem-finite-1}, there exists a $\mathbf{P}\in\mathcal{D}\left(\mathbb{S}^{n}\right)$ such that, for any $(x,i)\in\mathbb{R}^{n}\times \mathcal{S}$, $V^{0}(x,i)=\big<P(i)x,x\big>$.
Applying the dynamic programming principle to Problem (M-LQ)$^{0}$ and using It\^o's formula to $\left<P(\alpha_{t})X^{0}(t),X^{0}(t)\right>$, one has
\begin{align*}
&\big<P(i)x,x\big>
\leq \mathbb{E}\!\Big\{\!
    \int_{0}^{t}\!
    \left<\!
    \left(\!
    \begin{matrix}
    Q(\alpha)&S(\alpha)^{\top}\\
    S(\alpha)&R(\alpha)
    \end{matrix}
    \!\right)\!
    \left(\!
    \begin{matrix}
    X^{0}\\
    u
    \end{matrix}
    \!\right)\!,
    \!\left(\!
    \begin{matrix}
    X^{0}\\
    u
    \end{matrix}
    \!\right)
    \!\right>ds\\
&\quad+\big<P(\alpha_{t})X^{0}(t),X^{0}(t)\big>\Big\}\\
&= \mathbb{E}\!
    \int_{0}^{t}\!
    \left<\!
    \left(\!
    \begin{matrix}
    \mathcal{M}(P,\alpha)&\mathcal{L}(P,\alpha)\\
    \mathcal{L}(P,\alpha)^{\top}&R(\alpha)
    \end{matrix}
     \!\right) \!
    \left( \!
    \begin{matrix}
    X^{0}\\
    u
    \end{matrix}
     \!\right) \!,
     \!\left( \!
    \begin{matrix}
    X^{0}\\
    u
    \end{matrix}
    \!\right)
    \!\right>ds\\
    &\quad +\big<P(i)x,x\big>,
\end{align*}
which implies, for any $(x,i,u)\in\mathbb{R}^{n}\times\mathcal{S}\times\mathcal{U}_{ad}(x,i)$, the following holds:
\begin{align*}
\mathbb{E}\!\int_{0}^{t}\!
    \left<\!
    \left(\!
    \begin{matrix}
    \mathcal{M}(P,\alpha)&\mathcal{L}(P,\alpha) \\
    \mathcal{L}(P,\alpha)^{\top}&R(\alpha)
    \end{matrix}
    \!\right)\!
    \left(\!
    \begin{matrix}
    X^{0}\\
    u
    \end{matrix}
    \!\right)\!,
    \left(\!
    \begin{matrix}
    X^{0}\\
    u
    \end{matrix}
    \!\right)\!
    \right>ds\geq 0.
\end{align*}
Dividing both sides by $t$ and letting $t\downarrow 0$, we obtain $\mathbf{P}\in\mathcal{G}$.
On the other hand, for any $\mathbf{\widetilde{P}}\in\mathcal{G}$, we have
$$\big<P(i)x,x\big>=V^{0}(x,i)\geq \big<\widetilde{P}(i)x,x\big>, \qquad \forall (x,i)\in\mathbb{R}^{n}\times \mathcal{S}. $$
Therefore,
$$P(i)\geq\widetilde{P}(i),\qquad \forall i\in\mathcal{S}.$$
Consequently, $\mathbf{P}$ is the maximal element of $\mathcal{G}$. This completes the proof.

The following theorem provides a characterization for the closed-loop solvability of Problem (M-LQ).
\begin{thm}\label{thm-SLQ_0-closed}
The Problem (M-LQ) is closed-loop solvable if and only if the following conditions hold:
\begin{description}
  \item[(i)] The CAREs:
\begin{equation}\label{CAREs-LQ}
    \mathcal{M}(P,i)-\mathcal{L}(P,i) R(i)^{-1} \mathcal{L}(P,i)^{\top} = 0,
\end{equation}
admits a  solution $\mathbf{P}\in\mathcal{D}\left(\mathbb{S}^{n}\right)$ such that
\begin{align*}
  &\left[\!-R(\!1\!)^{-1}\!\mathcal{L}(\!P\!,\!1\!)^{\top}\!,\!-R(\!2\!)^{-1}\!\mathcal{L}(\!P\!,\!2\!)^{\top}\!,...,\!-R(\!D\!)^{-1}\!\mathcal{L}(\!P\!,\!D\!)^{\top}\!\right]\\
 &\triangleq-\mathbf{R}^{-1}\mathcal{L}(\mathbf{P})^{\top}\in\mathcal{H}_{\alpha};
\end{align*}
\item[(ii)] The  BSDE:
  \begin{equation}\label{eta}
      \begin{aligned}
        d\eta&=-\big\{\big[A(\alpha)^{\top}-\mathcal{L}(P,\alpha)R(\alpha)^{-1}B(\alpha)^{\top}\big]\eta+P(\alpha)b\\
       &\quad -\mathcal{L}(P,\alpha)R(\alpha)^{-1}\rho+q\big\}dt
        +\mathbf{z}\cdot d\mathbf{\widetilde{N}}(t),
      \end{aligned}
  \end{equation}
   admits a solution  $\left(\eta,\mathbf{z}\right)\in L_{\mathbb{F}}^{2}(\mathbb{R}^{n})\times \mathcal{D}\left(L_{\mathcal{P}}^{2}(\mathbb{R}^{n})\right)$.
\end{description}
In this case, the closed-loop optimal control $(\widehat{\mathbf{\Theta}},\widehat{\nu})$ admits the following representation:
\begin{equation}\label{closed-loop-saddle-point-3}
  \left\{
  \begin{aligned}
    \widehat{\Theta}(i)&=-R(i)^{-1}\mathcal{L}(P,i)^{\top},\quad\forall i\in\mathcal{S},\\
    \widehat{\nu}&=-R(\alpha)^{-1}\widetilde{\rho},
  \end{aligned}
  \right.
\end{equation}
where $\widetilde{\rho}\triangleq B(\alpha)^{\top}\eta+\rho$. Moreover, the value function is given by
\begin{equation}\label{LQ-value-function}
    \begin{aligned}
      V(x,i)&=\big<P(i)x,x\big>+\mathbb{E}\big\{2\big<\eta(0),x\big>+\int_{0}^{\infty}
      \big[2\big<\eta,b\big>\\
      &\quad-\big<R(\alpha)^{-1}\widetilde{\rho},\widetilde{\rho}\big>]dt\big\}.
    \end{aligned}
\end{equation}
\end{thm}

\textbf{Proof.}
\textbf{Necessity.} Let $(\widehat{\mathbf{\Theta}},\widehat{\nu})$ be a closed-loop optimal control. Then, by Corollary \ref{coro-SLQ-0}, $(\widehat{\mathbf{\Theta}},0)$ is a closed-loop optimal control for Problem (M-LQ)$^{0}$. Consequently, the Problem (M-LQ)$^{0}$ is well-posedness and by Lemma \ref{lem-finite-1}, there exists a maximal element $\mathbf{P}\in\mathcal{G}$ such that
\begin{equation}\label{SLQ_0-closed-p1}
    V^{0}(x,i)=\left<P(i)x,x\right>,\qquad\forall (x,i)\in\mathbb{R}^{n}\times\mathcal{S},
\end{equation}
and
$$
\left[
\begin{array}{cc}
  \mathcal{M}(P,i)  &  \mathcal{L}(P,i)  \\
  \mathcal{L}(P,i)^{\top}  & R(i)
\end{array}
\right]\geq 0,\quad\forall i\in\mathcal{S}.
$$
It follows from Schur's Lemma (see \cite{boyd1994linear}) that
\begin{equation}\label{SLQ_0-closed-p2}
 \mathcal{M}(P,i)-\mathcal{L}(P,i) R(i)^{-1} \mathcal{L}(P,i)^{\top} \geq 0, \quad\forall i\in\mathcal{S},
\end{equation}
Let  $X\equiv X_{\widehat{\Theta}}^{0}(\cdot;x,i,0)$. The Remark \ref{rem-useful} yields
\begin{align*}
&\big<P(i)x,x\big>=J^{0}(x,i;\widehat{\Theta}(\alpha)X)=\big<P(i)x,x\big> \\
&+\!\mathbb{E}\!\int_{0}^{\infty}\!
\big<\!\big[\mathcal{M}(P,\alpha)\!-\!\mathcal{L}(P,\alpha)R(\alpha)^{-1}\mathcal{L}(P,\alpha)^{\top}\big]X,X\!\big>dt \\
&+\mathbb{E}\!\int_{0}^{\infty}\!\big<\!R(\alpha)\!\big[\widehat{\Theta}(\alpha)
\!+\!R(\alpha)^{-1}\mathcal{L}(P,\alpha)^{\top}\!\big]X\!,\big[\widehat{\Theta}(\alpha)\\
&+\!R(\alpha)^{-1}\mathcal{L}(P,\alpha)^{\top}\!\big]X\!\big>dt.
\end{align*}
Combining with \eqref{SLQ_0-closed-p2}, we have
\begin{equation}\label{p1}
  \begin{aligned}
&\mathcal{M}(P,\alpha_{t})-\mathcal{L}(P,\alpha_{t}) R(\alpha_{t})^{-1} \mathcal{L}(P,\alpha_{t})^{\top}=0,\\
& \widehat{\Theta}(\alpha)=-R(\alpha)^{-1}\mathcal{L}(P,\alpha)^{\top},
\end{aligned}
\end{equation}
which implies $\mathbf{P}$ solves the CAREs \eqref{CAREs-LQ} and $-\mathbf{R}^{-1}\mathcal{L}(\mathbf{P})^{\top}\in\mathcal{H}_{\alpha}$.
On the other hand, let
{\scriptsize
$$
\left(\widehat{X},\widehat{Y},\widehat{\mathbf{\Gamma}}\right)\equiv
\left(X_{\widehat{\Theta}}(\cdot;x,i,\widehat{\nu}),Y_{\widehat{\Theta}}(\cdot;x,i,\widehat{\nu}),\mathbf{\Gamma}_{\widehat{\Theta}}(\cdot;x,i,\widehat{\nu})\right),
$$
}
and
\begin{equation}\label{eta-2}
\left\{
\begin{aligned}
&\eta(t)=\widehat{Y}(t)-P(\alpha_t)\widehat{X}(t),\\
&z_{j}(t)=\widehat{\Gamma}_{j}(t)-\big[P(j)-P(\alpha_{t-})\big]\widehat{X}(t),\text{ } t\geq 0.
\end{aligned}
\right.
\end{equation}
Then, $\left(\eta,\mathbf{z}\right)\in L_{\mathbb{F}}^{2}(\mathbb{R}^{n})\times\mathcal{D}\left(L_{\mathcal{P}}^{2}(\mathbb{R}^{n})\right)$. Noting stationary condition \eqref{SLQ-stationarity}, we have
\begin{align*}
 &0=B(\alpha)^{\top}\widehat{Y}+S_{\widehat{\Theta}}(\alpha)\widehat{X}+R(\alpha)\widehat{\nu}+\rho\\
 &=\!B(\alpha)^{\top}\!\big[\eta\!+\!P(\alpha)\widehat{X}\big]\!+\!\big[S(\alpha)\!+\!R(\alpha)\widehat{\Theta}(\alpha)\big]\!\widehat{X}\!+\!R(\alpha)\widehat{\nu}\!+\!\rho\\
 &=\big[\mathcal{L}(P,\alpha)^{\top}+R(\alpha)\widehat{\Theta}(\alpha)\big]\widehat{X}+R(\alpha)\widehat{\nu}+\widetilde{\rho}\\
 &=R(\alpha)\widehat{\nu}+\widetilde{\rho},
\end{align*}
which implies
\begin{equation}\label{p2}
\widehat{\nu}=-R(\alpha)^{-1}\widetilde{\rho},
\end{equation}
Then one can easily verify that $(\eta,\mathbf{z})$ defined in \eqref{eta-2} solves BSDE \eqref{eta} based on the relations \eqref{p1} and \eqref{p2}.

\textbf{Sufficiency.} Let $(\widehat{\mathbf{\Theta}},\widehat{\nu})$ be given by \eqref{closed-loop-saddle-point-3} and $\mathbf{P}\in\mathcal{D}\left(\mathbb{S}^{n}\right)$ be the solution to \eqref{CAREs-LQ}. Then we have
{\scriptsize
\begin{align*}
 \mathcal{M}(P,\alpha)\!+\!\mathcal{L}(P,\alpha)\widehat{\Theta}(\alpha)\!+\!\widehat{\Theta}(\alpha)^{\top}\mathcal{L}(P,\alpha)^{\top}\!+\!\widehat{\Theta}(\alpha)^{\top}R(\alpha)\widehat{\Theta}(\alpha)=0.
\end{align*}
}%
For any $\nu\in L_{\mathcal{P}}^{2}(\mathbb{R}^{m})$, set $X\equiv X_{\widehat{\Theta}}(\cdot;x,i,\nu)$.  Similar to Lemma \ref{lem-useful},  one has
\begin{equation}\label{main-result-1}
\begin{aligned}
&J_{\widehat{\Theta}}(x,i;\nu)=\mathbb{E}\int_{0}^{\infty}\big[2\big<P(\alpha)b+q+\widehat{\Theta}(\alpha)^{\top}\rho,X\big>\\
&\quad+\big<R(\alpha)\nu+2\rho,\nu\big>\big]dt+\big<P(i)x,x\big>.
\end{aligned}
\end{equation}
On the other hand, applying the It\^o's rule to $\big<\eta,X\big>$ and by some straightforward calculations, we have
\begin{align*}
 \mathbb{E}\big[\big<\eta(0),x\big>\big]&=\mathbb{E}\int_{0}^{\infty}\big\{
\big<P(\alpha)b+q+\widehat{\Theta}(\alpha)^{\top}\rho,X\big>\\
&\quad-\big<B(\alpha)\nu+b,\eta\big>\big\}dt.
\end{align*}
Combining with \eqref{main-result-1}, one can obtain
\begin{align*}
 &\quad J_{\widehat{\Theta}}(x,i;\nu)-\big<P(i)x,x\big>-2\mathbb{E}\big[\big<\eta(0),x\big>\big]\\
 &=\mathbb{E}\int_{0}^{\infty}\big[
\big<R(\alpha)\nu,\nu\big>+2\big<\widetilde{\rho},\nu\big>+2\big<b,\eta\big>\big]dt\\
 &=\mathbb{E}\int_{0}^{\infty}\big[
\big<R(\alpha)\nu,\nu\big>-2\big<R(\alpha)\widehat{\nu},\nu\big>+2\big<b,\eta\big>\big]dt\\
 &=\mathbb{E}\int_{0}^{\infty}\!\big[\!
\big<\!R(\alpha)(\nu\!-\!\widehat{\nu}),\nu\!-\!\widehat{\nu}\!\big>\!-\!\big<\!R(\alpha)\widehat{\nu},\widehat{\nu}\!\big>\!+\!2\big<b,\eta\big>\big]dt.
\end{align*}
Consequently, noting that $R(\alpha)>0$, we have
\begin{align*}
 &\quad J_{\widehat{\Theta}}(x,i;\nu)-J_{\widehat{\Theta}}(x,i;\widehat{\nu})\\
 &=\mathbb{E}\int_{0}^{\infty}\big<R(\alpha)\big(\nu-\widehat{\nu}\big),\nu-\widehat{\nu}\big> dt \geq 0,
\end{align*}
for all $(x,i,\nu)\in\mathbb{R}^{n}\times \mathcal{S}\times L_{\mathcal{P}}^{2}(\mathbb{R}^{m})$.
Therefore, by Proposition \ref{prop-LQ-control-characterization}, $(\mathbf{\widehat{\Theta}},\widehat{\nu})$ defined by \eqref{closed-loop-saddle-point-3} is a closed-loop optimal control of Problem (M-LQ). In addition, note that
\begin{align*}
\big<R(\alpha)\widehat{\nu},\widehat{\nu}\big> =\big<R(\alpha)^{-1}R(\alpha)\widehat{\nu},R(\alpha)\widehat{\nu}\big>
=\big<R(\alpha)^{-1}\widetilde{\rho},\widetilde{\rho}\big>.
\end{align*}
Hence, the value function of Problem (M-LQ) is given by
\begin{align*}
    V(x,i)&=J_{\widehat{\Theta}}(x,i;\widehat{\nu})=\big<P(i)x,x\big>+2\mathbb{E}\big[\big<\eta(0),x\big>\big]\\
   &\quad +\mathbb{E}\int_{0}^{\infty}\big\{2\big<b,\eta\big>-\big<R(\alpha)^{-1}\widetilde{\rho},\widetilde{\rho}\big>\big\}dt.
\end{align*}
This completes the proof.


The next corollary followed from Theorem \ref{thm-SLQ_0-closed} by setting $b,\,\sigma,\,q,\,\rho=0$.

\begin{cor}\label{coro-LQ-main-result-0}
The Problem (M-LQ)$^{0}$ is closed-loop solvable if and only if the CAREs \eqref{CAREs-LQ}
admits a  solution $\mathbf{P}\in\mathcal{D}\left(\mathbb{S}^{n}\right)$ such that
$-\mathbf{R}^{-1}\mathcal{L}(\mathbf{P})^{\top}\in\mathcal{H}_{\alpha}.$
In this case, the closed-loop optimal control $(\widehat{\mathbf{\Theta}},\widehat{\nu})$ admits the following representation:
\begin{equation}\label{closed-loop-saddle-point-4}
  \left\{
  \begin{aligned}
    &\widehat{\Theta}(i)=-R(i)^{-1}\mathcal{L}(P,i)^{\top},\quad\forall i\in\mathcal{S},\\
    &\widehat{\nu}=0.
  \end{aligned}
  \right.
\end{equation}
Moreover, the value function is given by
\begin{equation}\label{LQ-value-function-0}
      V(x,i)=\big<P(i)x,x\big>.
\end{equation}
\end{cor}

\section{The closed-loop solvable of Problem (M-GLQ)}\label{section-GLQ}
In this section, we return to study the Nash equilibrium game problem introduced in Section \ref{section-introduction}.
For any $u_{2}\in L_{\mathcal{P}}^{2}(\mathbb{R}^{m_{2}})$ and $u_{1}\in L_{\mathcal{P}}^{2}(\mathbb{R}^{m_{1}})$, we begin by introducing the following notations
$$\begin{aligned}
&\mathcal{U}_{ad}^{1}(x,i;u_{2})\triangleq\left\{u_{1}\in L_{\mathcal{P}}^{2}(\mathbb{R}^{m_{1}})\big|\left(u_{1},u_{2}\right)\in\mathcal{U}_{ad}(x,i)\right\},\\
&\mathcal{U}_{ad}^{2}(x,i;u_{1})\triangleq\left\{u_{2}\in L_{\mathcal{P}}^{2}(\mathbb{R}^{m_{2}})|\left(u_{1},u_{2}\right)\in\mathcal{U}_{ad}(x,i)\right\}.
\end{aligned}$$
For any given $\mathbf{\Theta}\equiv (\mathbf{\Theta_{1}},\mathbf{\Theta_{2}})\in\mathcal{D}\left(\mathbb{R}^{m_{1}\times n}\right)\times \mathcal{D}\left(\mathbb{R}^{m_{2}\times n}\right)$ and $(\nu_{1},\nu_{2})\in L_{\mathcal{P}}^{2}(\mathbb{R}^{m_{1}})\times L_{\mathcal{P}}^{2}(\mathbb{R}^{m_{2}})$, let
{\scriptsize
\begin{equation}\label{GLQ-state-closed}
   \left\{
   \begin{aligned}
   \dot{X}_{\Theta}(t)&\!=\!A_{\Theta}(\alpha_{t})X_{\Theta}(t)\!+\!B_{1}(\alpha_{t})\nu_{1}(t)\!+\!B_{2}(\alpha_{t})\nu_{2}(t)\!+\!b(t),\\
   X_{\Theta}(0)&=x,\quad \alpha(0)=i,\quad t\geq0,
   \end{aligned}
   \right.
 \end{equation}
 }%
 and for any $(\mathbf{\Theta}_{l},\nu_{l};u_{k})\in\mathcal{D}\left(\mathbb{R}^{m_{l}\times n}\right)\times L_{\mathcal{P}}^{2}(\mathbb{R}^{m_{l}})\times  L_{\mathcal{P}}^{2}(\mathbb{R}^{m_{k}})$, let
 \begin{equation}\label{GLQ-state-k}
   \left\{
   \begin{aligned}
   \dot{X}_{k}(t)&=A_{k}(\alpha_{t})X_{k}(t)+B_{k}(\alpha_{t})u_{k}(t)+b_{k}(t),\\
   X_{k}(0)&=x,\quad \alpha(0)=i,\quad t\geq0,\quad k=1,2,
   \end{aligned}
   \right.
 \end{equation}
 where $k,l\in\{1,2\}$, $k\neq l$ and
\begin{equation}
\left\{
\begin{aligned}
&A_{\Theta}(i)=A(i)+B_{1}(i)\Theta_{1}(i)+B_{2}(i)\Theta_{2}(i),\\
&A_{k}(i)=A(i)+B_{l}(i)\Theta_{l}(i),\\
&b_{k}=B_{l}(\alpha)\nu_{l}+b.\\
\end{aligned}
\right.
\end{equation}
We denote the solution of \eqref{GLQ-state-closed} as $X_{\Theta}(\cdot;x,i,\nu_1,\nu_2)$ and the solution of \eqref{GLQ-state-k} as $ X_{k}(\cdot;x,i,\mathbf{\Theta}_{l},\nu_{l},u_{k})$.
Similar to Definition \ref{def-stabilizable}, we introduce the following concept.
\begin{defn}\label{def-GLQ-stabilizable}
The system \eqref{GLQ-state} is said to be $L^{2}$-stabilizable if and only if there exists a pair $(\mathbf{\Theta_{1}},\mathbf{\Theta_{2}})\in\mathcal{D}\big(\mathbb{R}^{m_{1}\times n}\big)\times \mathcal{D}\big(\mathbb{R}^{m_{2}\times n}\big)$ such that
\begin{equation}\label{L-3}
\begin{aligned}
&P(i)[A(i)\!+\!B_{1}(i)\Theta_{1}(i)\!+\!B_{2}(i)\Theta_{2}(i)]\!+\![A(i)\!+\!B_{1}(i)\Theta_{1}(i)\\
&+\!B_{2}(i)\Theta_{2}(i)]^{\top}P(i)\!+\!\sum_{j=1}^{D}\pi_{ij}P(j)<0,\quad i\in\mathcal{S},
\end{aligned}
\end{equation}
for some $\mathbf{P}\in \mathcal{D}\left(\mathbb{S}_{+}^{n}\right)$.
Any pair $(\mathbf{\Theta_{1}},\mathbf{\Theta_{2}})\in\mathcal{D}\big(\mathbb{R}^{m_{1}\times n}\big)\times \mathcal{D}\big(\mathbb{R}^{m_{2}\times n}\big)$  satisfying the above condition is called a stabilizer of system \eqref{GLQ-state} and we denote the $\widehat{\mathcal{H}}_{\alpha}$ represents the set consisting of all stabilizers.
\end{defn}

Next, we will work under the following assumption.\\
\textbf{(H2)} System  \eqref{GLQ-state} is $L^{2}$-stabilizable, i.e, $\widehat{\mathcal{H}}_{\alpha}\neq\emptyset$.

\begin{defn}\label{def-closed-loop-equilibrium-point-GLQ}
A 4-tuple $(\mathbf{\widehat{\Theta}_{1}},\widehat{\nu}_{1};\mathbf{\widehat{\Theta}_{2}},\widehat{\nu}_{2})\in\mathcal{D}\left(\mathbb{R}^{m_{1}\times n}\right)\times L_{\mathcal{P}}^{2}(\mathbb{R}^{m_{1}})\times \mathcal{D}\left(\mathbb{R}^{m_{1}\times n}\right)\times L_{\mathcal{P}}^{2}(\mathbb{R}^{m_{2}})$ is called a closed-loop Nash equilibrium point of Problem (M-GLQ) if:\\
$\mathbf{(i)}$ $(\mathbf{\widehat{\Theta}_{1}},\mathbf{\widehat{\Theta}_{2}})\in \widehat{\mathcal{H}}_{\alpha}$,\\
$\mathbf{(ii)}$ for any $(x,i)\in\mathbb{R}^{n}\times\mathcal{S}$ and $(u_{1},u_{2})\in\mathcal{U}_{ad}^{1}(x,i;\widehat{\Theta}_{2}(\alpha)\widehat{X}_{1}+\widehat{\nu}_{2})\times \mathcal{U}_{ad}^{2}(x,i;\widehat{\Theta}_{1}(\alpha)\widehat{X}_{2}+\widehat{\nu}_{1})$, the following holds:
  \begin{equation}\label{closed-loop-Nash-equilibrium-point}
  \left\{
      \begin{aligned}
&J_{1}(x,i;\widehat{\Theta}_{1}(\alpha)\widehat{X}_{\Theta}+\widehat{\nu}_{1},\widehat{\Theta}_{2}(\alpha)\widehat{X}_{\Theta}+\widehat{\nu}_{2})\\
&\leq J_{1}(x,i;u_{1},\widehat{\Theta}_{2}(\alpha)\widehat{X}_{1}+\widehat{\nu}_{2}),\\
&J_{2}(x,i;\widehat{\Theta}_{1}(\alpha)\widehat{X}_{\Theta}+\widehat{\nu}_{1},\widehat{\Theta}_{2}(\alpha)\widehat{X}_{\Theta}+\widehat{\nu}_{2})\\
&\leq J_{2}(x,i;\widehat{\Theta}_{1}(\alpha)\widehat{X}_{2}+\widehat{\nu}_{1},u_{2}),
      \end{aligned}
      \right.
  \end{equation}
where
\begin{equation*}
\left\{
    \begin{aligned}
    &\widehat{X}_{\Theta}=X_{\widehat{\Theta}}(\cdot;x,i,\widehat{\nu}_{1},\widehat{\nu}_{2}),\\
    &\widehat{X}_{1}=X_{1}(\cdot;x,i,u_1,\mathbf{\widehat{\Theta}_{2}},\widehat{\nu}_{2}),\\
    &\widehat{X}_{2}=X_{2}(\cdot;x,i,\mathbf{\widehat{\Theta}_{1}},\widehat{\nu}_{1},u_2).
    \end{aligned}
    \right.
\end{equation*}
In addition, Problem (M-GLQ) is said to be (uniquely) closed-loop solvable if it has a (unique)  closed-loop Nash equilibrium point.
\end{defn}

Clearly, by some straightforward calculations, one has:
\begin{equation}\label{cost-GLQ}
    \begin{aligned}
        & J_{1}\left(x,i;u_{1},\widehat{\Theta}_{2}(\alpha)\widehat{X}_{1}+\widehat{\nu}_{2}\right)\\
        =&\widehat{J}_{1}\left(x,i;u_{1}\right)
        +\mathbb{E}\int_{0}^{\infty}\big<R_{22}^{1}(\alpha)\widehat{\nu}_{2}+2\rho_{2}^{1},\widehat{\nu}_{2}\big>dt,\\
        &J_{2}\left(x,i;\widehat{\Theta}_{1}(\alpha)\widehat{X}_{2}+\widehat{\nu}_{1},u_{2}\right)\\
        =&\widehat{J}_{2}\left(x,i;u_{2}\right)
        +\mathbb{E}\int_{0}^{\infty}\big<R_{11}^{2}(\alpha)\widehat{\nu}_{1}+2\rho_{1}^{2},\widehat{\nu}_{1}\big>dt,
    \end{aligned}
\end{equation}
where
{\scriptsize
\begin{equation}\label{cost-GLQ-1}
\begin{aligned}
    \widehat{J}_{1}\left(x,i;u_{1}\right)
    &=\mathbb{E}\int_{0}^{\infty}\left[
    \left<
    \left(
    \begin{matrix}
    \widehat{Q}^{1}(\alpha)&\widehat{S}^{1}(\alpha)^{\top}\\
    \widehat{S}^{1}(\alpha)&R_{11}^{1}(\alpha)
    \end{matrix}
    \right)
    \left(
    \begin{matrix}
    \widehat{X}_{1}\\
    u_{1}
    \end{matrix}
    \right),
    \left(
    \begin{matrix}
    \widehat{X}_{1}\\
    u_{1}
    \end{matrix}
    \right)
    \right>\right.\\
    &\quad\left.+2\left<
    \left(
    \begin{matrix}
    \widehat{q}^{1}\\
    \widehat{\rho}^{1}
    \end{matrix}
    \right),
    \left(
    \begin{matrix}
    \widehat{X}_{1}\\
    u_{1}
    \end{matrix}
    \right)
    \right>\right]dt,
  \end{aligned}
\end{equation}
}
{\scriptsize
\begin{equation}\label{cost-GLQ-2}
\begin{aligned}
    \widehat{J}_{2}\left(x,i;u_{2}\right)
    & = \mathbb{E}\int_{0}^{\infty}\left[
    \left<
    \left(
    \begin{matrix}
    \widehat{Q}^{2}(\alpha)&\widehat{S}^{2}(\alpha)^{\top}\\
    \widehat{S}^{2}(\alpha)&R_{22}^{2}(\alpha)
    \end{matrix}
    \right)
    \left(
    \begin{matrix}
    \widehat{X}_{2}\\
    u_{2}
    \end{matrix}
    \right),
    \left(
    \begin{matrix}
    \widehat{X}_{2}\\
    u_{2}
    \end{matrix}
    \right)
    \right>\right.\\
    &\quad\left.+2\left<
    \left(
    \begin{matrix}
    \widehat{q}^{2}\\
    \widehat{\rho}^{2}
    \end{matrix}
    \right),
    \left(
    \begin{matrix}
    \widehat{X}_{2}\\
    u_{2}
    \end{matrix}
    \right)
    \right>\right]dt,
  \end{aligned}
\end{equation}
}
and
{\scriptsize
\begin{equation}\label{notation-GLQ}
\begin{aligned}
&\widehat{Q}^{1}(i)\!= \!Q^{1}(i)\!+\!S_{2}^{1}(i)^{\top}\widehat{\Theta}_{2}(i)
\!+\!\widehat{\Theta}_{2}(i)^{\top}S_{2}^{1}(i)
\!+\!\widehat{\Theta}_{2}(i)^{\top}R_{22}^{1}(i)\widehat{\Theta}_{2}(i),\\
&\widehat{S}^{1}(i)\!=\!S_{1}^{1}(i)\!+\!R_{12}^{1}(i)\widehat{\Theta}_{2}(i),\quad
\widehat{\rho}^{1}\!=\!\rho_{1}^{1}\!+\!R_{12}^{1}(\alpha)\widehat{\nu}_{2},\\
&\widehat{q}^{1}\!=\!q^{1}\!+\!\widehat{\Theta}_{2}(\alpha)^{\top}\rho_{2}^{1}
\!+\!\big[S_{2}^{1}(\alpha)^{\top}\!+\!\widehat{\Theta}_{2}(\alpha)^{\top}R_{22}^{1}(\alpha)\big]\widehat{\nu}_{2},\\
&\widehat{Q}^{2}(i)\!= \!Q^{2}(i)\!+\!S_{1}^{2}(i)^{\top}\widehat{\Theta}_{1}(i)
\!+\!\widehat{\Theta}_{1}(i)^{\top}S_{1}^{2}(i)
\!+\!\widehat{\Theta}_{1}(i)^{\top}R_{11}^{2}(i)\widehat{\Theta}_{1}(i),\\
&\widehat{S}^{2}(i)\!=\!S_{2}^{2}(i)\!+\!R_{21}^{2}(i)\widehat{\Theta}_{1}(i),\quad
\widehat{\rho}^{2}\!=\!\rho_{2}^{2}\!+\!R_{21}^{2}(\alpha)\widehat{\nu}_{1},\\
&\widehat{q}^{2}\!=\!q^{2}\!+\!\widehat{\Theta}_{1}(\alpha)^{\top}\rho_{1}^{2}
\!+\!\big[S_{1}^{2}(\alpha)^{\top}\!+\!\widehat{\Theta}_{1}(\alpha)^{\top}R_{11}^{2}(\alpha)\big]\widehat{\nu}_{1}.
\end{aligned}
\end{equation}
}

Therefore, combining with the Definition \ref{def-closed-loop-equilibrium-point-GLQ}, a 4-tuple $(\mathbf{\widehat{\Theta}_{1}},\widehat{\nu}_{1};\mathbf{\widehat{\Theta}_{2}},\widehat{\nu}_{2})$ is a closed-loop Nash equilibrium point of Problem (M-GLQ) if and only if the following holds:
\begin{description}
  \item[(i)] $(\mathbf{\widehat{\Theta}_{1}},\mathbf{\widehat{\Theta}_{2}})\in \widehat{\mathcal{H}}_{\alpha}$,
  \item[(ii)] $(\mathbf{\widehat{\Theta}_{1}},\widehat{\nu}_{1})$ (respectively, $(\mathbf{\widehat{\Theta}_{2}},\widehat{\nu}_{2})$) is a closed-loop optimal control for problem with cost functional $\widehat{J}_{1}\left(x,i;u_{1}\right)$ (respectively, $\widehat{J}_{2}\left(x,i;u_{2}\right)$) and state constraint $\widehat{X}_{1}$ (respectively, $\widehat{X}_{2}$).
\end{description}

 To simplify our further analysis, for any given $(\mathbf{P_{1}},\mathbf{P_{2}})\in\mathcal{D}\left(\mathbb{S}^{n}\right)\times \mathcal{D}\left(\mathbb{S}^{n}\right)$, we introduce the following notations:
 $$
\begin{array}{l}
\widehat{A}_{k}(i)\!\triangleq \!A(i)\!+\!B_{l}(i)\widehat{\Theta}_{l}(i),\quad \widehat{b}_{k}\!\triangleq \!B_{l}(\alpha)\widehat{\nu}_{l}\!+\!b\\
\widehat{\mathcal{L}}_{k}(P_{k},i)\!\triangleq\! P_{k}(i)B_{k}(i)\!+\!\widehat{S}^{k}(i)^{\top},\\
\widehat{\mathcal{M}}_{k}(P_{k},i)\!\triangleq \!P_{k}(i)\widehat{A}_{k}(i)\!+\!\widehat{A}_{k}(i)^{\top}P_{k}(i)\!+\!\widehat{Q}^{k}(i)\\
\qquad\qquad\quad+\!\sum_{j=1}^{D}\!\pi_{ij}P_{k}(j),
\text{ } k,l\in\{1,2\}, \text{ } l\neq k.
\end{array}
$$
Then
\begin{equation}\label{relation}
\left\{
\begin{aligned}
&\widehat{\mathcal{M}}_{1}(P_{1},i)=\mathcal{M}_{1}(P_{1},i)+\mathcal{L}_2^1(P_{1},i)\widehat{\Theta}_2(i)\\
&\qquad+\widehat{\Theta}_2(i)^{\top}\mathcal{L}_2^1(P_{1},i)^{\top}+\widehat{\Theta}_2(i)^{\top}R_{22}^{1}(i)\widehat{\Theta}_2(i),\\
&\widehat{\mathcal{M}}_{2}(P_{2},i)=\mathcal{M}_{2}(P_{2},i)+\mathcal{L}_1^2(P_{2},i)\widehat{\Theta}_1(i)\\
&\qquad+\widehat{\Theta}_1(i)^{\top}\mathcal{L}_1^2(P_{2},i)^{\top}+\widehat{\Theta}_1(i)^{\top}R_{11}^{2}(i)\widehat{\Theta}_1(i),\\
&\widehat{\mathcal{L}}_{1}(P_{1},i)=\mathcal{L}_{1}^{1}(P_{1},i)+\widehat{\Theta}_{2}(i)^{\top}R_{21}^{1}(i),\\
&\widehat{\mathcal{L}}_{2}(P_{2},i)=\mathcal{L}_{2}^{2}(P_{2},i)+\widehat{\Theta}_{1}(i)^{\top}R_{12}^{2}(i).
\end{aligned}
\right.
\end{equation}
where
$$
\left\{
\begin{array}{l}
\mathcal{M}_{k}(P_{k},i)\triangleq P_{k}(i)A(i)+A(i)^{\top}P_{k}(i)\\
\quad\qquad\qquad+Q^{k}(i)+\sum_{j=1}^{D}\pi_{ij}P_{k}(j),\\
\mathcal{L}_{k}^{l}(P_{l},i)\triangleq P_{l}(i)B_{k}(i)+S_{k}^{l}(i)^{\top},\quad k,l=1,2.
\end{array}
\right.
$$

Consequently, by Theorem \ref{thm-SLQ_0-closed} and relation \eqref{relation},  we obtain the main result of this section.

\begin{thm}\label{thm-GLQ-closed}
 A 4-tuple $(\mathbf{\widehat{\Theta}_{1}},\!\widehat{\nu}_{1};\!\mathbf{\widehat{\Theta}_{2}},\!\widehat{\nu}_{2})\!\in\mathcal{D}\!(\mathbb{R}^{m_{1}\times n})$ $\times L_{\mathcal{P}}^{2}(\mathbb{R}^{m_{1}})\times \mathcal{D}\left(\mathbb{R}^{m_{1}\times n}\right)\times L_{\mathcal{P}}^{2}(\mathbb{R}^{m_{2}})$ is a closed-loop Nash equilibrium point of Problem (M-GLQ) if and only if:\\
 $\mathbf{(i)}$ $(\mathbf{\widehat{\Theta}_{1}},\mathbf{\widehat{\Theta}_{2}})\in \widehat{\mathcal{H}}_{\alpha}$,\\
 $\mathbf{(ii)}$ the following two cross-coupled CAREs:
  \begin{equation}\label{CAREs-GLQ-1}
  \begin{aligned}
  0&=\mathcal{M}_{1}(P_{1},i)-\mathcal{L}_{1}^{1}(P_{1},i)R_{11}^{1}(i)^{-1}\mathcal{L}_{1}^{1}(P_{1},i)^{\top}\\
  &+\big[\mathcal{L}_{2}^{1}(P_{1},i)-\mathcal{L}_{1}^{1}(P_{1},i)R_{11}^{1}(i)^{-1}R_{12}^{1}(i)\big]\widehat{\Theta}_{2}(i)\\
  &+\widehat{\Theta}_{2}(i)^{\top}\big[\mathcal{L}_{2}^{1}(P_{1},i)-\mathcal{L}_{1}^{1}(P_{1},i)R_{11}^{1}(i)^{-1}R_{12}^{1}(i)\big]^{\top}\\
  &+\widehat{\Theta}_{2}(i)^{\top}\big[R_{22}^{1}(i)-R_{21}^{1}(i)R_{11}^{1}(i)^{-1}R_{12}^{1}(i)\big]\widehat{\Theta}_{2}(i),
  \end{aligned}
  \end{equation}
   \begin{equation}\label{CAREs-GLQ-2}
    \begin{aligned}
  0&=\mathcal{M}_{2}(P_{2},i)-\mathcal{L}_{2}^{2}(P_{2},i)R_{22}^{2}(i)^{-1}\mathcal{L}_{2}^{2}(P_{2},i)^{\top}\\
  &+\big[\mathcal{L}_{1}^{2}(P_{2},i)-\mathcal{L}_{2}^{2}(P_{2},i)R_{22}^{2}(i)^{-1}R_{21}^{2}(i)\big]\widehat{\Theta}_{1}(i)\\
  &+\widehat{\Theta}_{1}(i)^{\top}\big[\mathcal{L}_{1}^{2}(P_{2},i)-\mathcal{L}_{2}^{2}(P_{2},i)R_{22}^{2}(i)^{-1}R_{21}^{2}(i)\big]^{\top}\\
  &+\widehat{\Theta}_{1}(i)^{\top}\big[R_{11}^{2}(i)-R_{12}^{2}(i)R_{22}^{2}(i)^{-1}R_{21}^{2}(i)\big]\widehat{\Theta}_{1}(i),
  \end{aligned}
  \end{equation}
  admit a solution $(\mathbf{P_{1}},\mathbf{P_{2}})\in\mathcal{D}\left(\mathbb{S}^{n}\right)\times \mathcal{D}\left(\mathbb{S}^{n}\right)$ such that
  \begin{equation}\label{CAREs-GLQ-constraint}
  \left\{
      \begin{aligned}
      &R_{11}^{1}(i)\widehat{\Theta}_{1}(i)+R_{12}^{1}(i)\widehat{\Theta}_{2}(i)+\mathcal{L}_{1}^{1}(P_{1},i)^{\top}=0,\\
      &R_{21}^{2}(i)\widehat{\Theta}_{1}(i)+R_{22}^{2}(i)\widehat{\Theta}_{2}(i)+\mathcal{L}_{2}^{2}(P_{2},i)^{\top}=0,
      \end{aligned}
      \right.
  \end{equation}
  \item[(iii)] the following two cross-coupled BSDEs:
{\scriptsize
\begin{equation}\label{GLQ-eta}
   \begin{aligned}
&d\eta_{1}=\big\{\big[\big(\mathcal{L}_{1}^{1}(P_{1},\alpha)+\widehat{\Theta}_{2}(\alpha)^{\top}R_{21}^{1}(\alpha)\big)R_{11}^{1}(\alpha)^{-1}B_{1}(\alpha)^{\top}\\
&\quad-\widehat{A}_{1}(\alpha)^{\top}\big]\eta_{1}-\big[S_{2}^{1}(\alpha)^{\top}+\widehat{\Theta}_{2}(\alpha)^{\top}R_{22}^{1}(\alpha)\big]\widehat{\nu}_{2}\\
&\quad +\big(\mathcal{L}_{1}^{1}(P_{1},\alpha)-\widehat{\Theta}_{2}(\alpha)^{\top}R_{21}^{1}(\alpha)\big)R_{11}^{1}(\alpha)^{-1}\big(\rho_{1}^{1}-R_{12}^{1}(\alpha)\widehat{\nu}_{2}\big)\\
&\quad-P_{1}(\alpha)\widehat{b}_{1}-q^{1}-\widehat{\Theta}_{2}(\alpha)^{\top}\rho_{2}^{1}\big\}dt+\mathbf{z}_{1}\cdot d\mathbf{\widetilde{N}}(t),\\
&d\eta_{2}=\big\{\big[\big(\mathcal{L}_{2}^{2}(P_{2},\alpha)+\widehat{\Theta}_{1}(\alpha)^{\top}R_{12}^{2}(\alpha)\big)R_{22}^{2}(\alpha)^{-1}B_{2}(\alpha)^{\top}\\
&\quad-\widehat{A}_{2}(\alpha)^{\top}\big]\eta_{2}-\big[S_{1}^{2}(\alpha)^{\top}+\widehat{\Theta}_{1}(\alpha)^{\top}R_{11}^{2}(\alpha)\big]\widehat{\nu}_{1}\\
&\quad+\big(\mathcal{L}_{2}^{2}(P_{2},\alpha)+\widehat{\Theta}_{1}(\alpha)^{\top}R_{12}^{2}(\alpha)\big)R_{22}^{2}(\alpha)^{-1}\big(\rho_{2}^{2}+R_{21}^{2}(\alpha)\widehat{\nu}_{1}\big)\\
&\quad-P_{2}(\alpha)\widehat{b}_{2}-q^{2}-\widehat{\Theta}_{1}(\alpha)^{\top}\rho_{1}^{2}\big\}dt+\mathbf{z}_{2}\cdot d\mathbf{\widetilde{N}}(t),
   \end{aligned}
  \end{equation}
}
admits a solution $(\eta_{1},\mathbf{z}_{1};\eta_{2},\mathbf{z}_{2})\in L_{\mathbb{F}}^{2}(\mathbb{R}^{n})\times\mathcal{D}\left(L_{\mathcal{P}}^{2}(\mathbb{R}^{n})\right)\times L_{\mathbb{F}}^{2}(\mathbb{R}^{n})\times\mathcal{D}\left(L_{\mathcal{P}}^{2}(\mathbb{R}^{n})\right)$ such that
\begin{equation}\label{CAREs-GLQ-eta-constraint}
  \left\{
      \begin{aligned}
      &R_{11}^{1}(\alpha)\widehat{\nu}_{1}+R_{12}^{1}(\alpha)\widehat{\nu}_{2}+B_1(\alpha)^{\top}\eta_{1}+\rho_{1}^{1}=0,\\
      &R_{21}^{2}(\alpha)\widehat{\nu}_{1}+R_{22}^{2}(\alpha)\widehat{\nu}_{2}+B_2(\alpha)^{\top}\eta_{2}+\rho_{2}^{2}=0.
      \end{aligned}
      \right.
  \end{equation}
\end{thm}

The following result is an obvious corollary of Theorem \ref{thm-GLQ-closed}.
\begin{cor}\label{coro-GLQ-0-closed}
If $(\mathbf{\widehat{\Theta}_{1}},\widehat{\nu}_{1};\mathbf{\widehat{\Theta}_{2}},\widehat{\nu}_{2})$ is a closed-loop Nash equilibrium point of Problem (M-GLQ), then $(\mathbf{\widehat{\Theta}_{1}},0;\mathbf{\widehat{\Theta}_{2}},0)$ is a closed-loop Nash equilibrium point of Problem (M-GLQ)$^{0}$.
\end{cor}

\section{Examples}\label{section-Examples}
In this section, we present two examples illustrating the results established in the previous sections. The numerical algorithms we have used are similar to those in \cite{Jianhui-Huang-2015,Li-Zhou-Rami-2003-ID-MLQ-IF}. Without loss of generality, we suppose the state space of the Markov chain $\alpha$ is $\mathcal{S}:=\left\{1,2,3\right\}$ and the generator is given by
$$\bf{\pi}=\left[\begin{array}{ccc}
     \pi_{11}   &   \pi_{12}   &   \pi_{13}\\
    \pi_{21}   &   \pi_{22}   &   \pi_{23}\\
    \pi_{31}   &   \pi_{32}   &   \pi_{33}
\end{array}\right]
=\left[\begin{array}{ccc}
     -0.5   &   0.2   &   0.3\\
    0.4   &   -0.6   &   0.2\\
    0.1   &   0.3   &   -0.4
\end{array}\right]$$
\begin{exmp}\label{exam-LQ}
Consider the following three-dimensional state equation
\begin{equation}\label{state-LQ-exam}
  \left\{
 \begin{aligned}
   \dot{X}(t)&=A\left(\alpha_{t}\right)X(t)+B\left(\alpha_{t}\right)u(t),\quad t\geq 0, \\
    X(0)&=x,\quad \alpha(0)=i,
   \end{aligned}
  \right.
\end{equation}
with the cost functional
\begin{equation}\label{cost-LQ-exam}
\begin{aligned}
    J(x,i;u)\!
    &\triangleq\! \mathbb{E}\!\int_{0}^{\infty}\!
    \left<\!
    \left(\!
    \begin{matrix}
    Q(\alpha)&S(\alpha)^{\top} \\
    S(\alpha)&R(\alpha)
    \end{matrix}
    \!\right)\!
    \left(\!
    \begin{matrix}
    X\\
    u
    \end{matrix}
    \!\right)\!,
    \left(\!
    \begin{matrix}
    X\\
    u
    \end{matrix}
    \!\right)
    \!\right>dt.
  \end{aligned}
\end{equation}
Let
\begin{align*}
 &A(1)=\left[\begin{matrix}
     -1   &   0   &   -1\\
    0   &   -2   &   0\\
    0   &   0   &   -1
\end{matrix}\right],\quad
A(2)=\left[\begin{matrix}
     -1   &   -1   &   0\\
    0   &   -1   &   -1\\
    0   &   0   &   -1
\end{matrix}\right],\\
&A(3)=\left[\begin{matrix}
     -2   &   0   &   0\\
    0   &   -1   &   0\\
    -1   &   0   &   -1
\end{matrix}\right],\quad
B(1)=\left[\begin{matrix}
     2   &   1   \\
    1   &   -2   \\
    0   &   3
\end{matrix}\right],\\
&B(2)=\left[\begin{matrix}
     1   &   -1   \\
    0   &   3   \\
    2   &   0
\end{matrix}\right],\quad
B(3)=\left[\begin{matrix}
     -2   &   2  \\
    0   &   4   \\
    1   &   0
\end{matrix}\right],\\
&\left[\begin{array}{cc}
    Q(1) & S(1) ^{\top}\\
    S(1) & R(1)
    \end{array}\right]
=\left[\begin{array}{ccccc}
     1.13   &   0.29   &   0.37   &   0.13   &   0.26\\
     0.29   &   0.15   &   0.29   &   0.13   &   0.07\\
     0.37   &   0.29   &   1.30   &   0.27   &   0.27\\
     0.13   &   0.13   &   0.27   &   0.23   &   0.11\\
     0.26   &   0.07   &   0.27   &   0.11   &   0.18
\end{array}\right],\\
&\left[\begin{array}{cc}
    Q(2) & S(2)^{\top} \\
    S(2) & R(2)
    \end{array}\right]
=\left[\begin{array}{ccccc}
     0.16   &   0.24   &   0.26   &   0.17   &   0.12\\
     0.24   &   0.53   &   0.71   &   0.26   &   0.29\\
     0.26   &   0.71   &   1.21   &   0.23   &   0.29\\
     0.17   &   0.26   &   0.23   &   0.24   &   0.17\\
     0.12   &   0.29   &   0.29   &   0.17   &   0.23
\end{array}\right],\\
&\left[\begin{array}{cc}
    Q(3) & S(3) ^{\top}\\
    S(3) & R(3)
    \end{array}\right]
=\left[\begin{array}{ccccc}
     0.39   &   0.27   &   0.24   &   0.24   &   0.22\\
     0.27   &   0.92   &   0.29   &   0.47   &   0.21\\
     0.24   &   0.29   &   0.39   &   0.27   &   0.21\\
     0.24   &   0.47   &   0.27   &   0.30   &   0.17\\
     0.22   &   0.21   &   0.21   &   0.17   &   0.19
\end{array}\right].
\end{align*}
Solving the corresponding  CAREs \eqref{CAREs-LQ}, we obtain
\begin{align*}
    P(1)&=\left[\begin{matrix}
     0.14589036  &  0.00191714  &  -0.05966332\\
     0.00191714  &  0.01377722  &   0.02024132\\
    -0.05966332  &  0.02024132  &   0.09509278
    \end{matrix}\right],\\
     P(2)&=\left[\begin{matrix}
     0.03965578  &  0.00277610  &  -0.02186540\\
     0.00277610  &  0.01259710  &   0.02536500\\
    -0.02186540  &  0.02536500  &   0.12162428
    \end{matrix}\right],\\
     P(3)&=\left[\begin{matrix}
     0.03530583  & -0.02129959  &   0.01944559\\
    -0.02129959  &  0.04517101  &  -0.03842815\\
     0.01944559  & -0.03842815  &   0.05171899
    \end{matrix}\right],
\end{align*}
with residuals
$$
\begin{aligned}
&\|E_{1}(\mathbf{P})\|=7.3846\times 10^{-8},\text{ }
\|E_{2}(\mathbf{P})\|= 4.6331\times 10^{-8},\\
&\|E_{3}(\mathbf{P})\|=8.8894\times 10^{-8}.
\end{aligned}$$
Here, the residual $E_{i}(\mathbf{P})$ is defined as follows:
\begin{equation}
   E_{i}(\mathbf{P})\triangleq \mathcal{M}(P,i)-\mathcal{L}(P,i) R(i)^{-1} \mathcal{L}(P,i)^{\top},\text{ } i=1,2,3.
\end{equation}
By Corollary \ref{coro-LQ-main-result-0}, the closed-loop optimal control $\mathbf{\widehat{\Theta}}=-\mathbf{R}^{-1}\mathcal{L}(\mathbf{P})^{\top}$ is given by:
\begin{align*}
    &\widehat{\Theta}(1)=\left[\begin{matrix}
    -1.76547241    &   -0.51230503   &    0.65870032\\
    -0.16035642   &   -0.27073948   &   -2.93105231
    \end{matrix}\right],\\
     &\widehat{\Theta}(2)=\left[\begin{matrix}
    -0.87788658   &   -0.64082962   &    -1.43972847\\
     0.26333999   &   -0.93945289   &    -0.62263722
    \end{matrix}\right],\\
    &\widehat{\Theta}(3)=\left[\begin{matrix}
    -0.03409914   &  -1.10028717    &    -1.33654538\\
    -1.05061288   &  -0.84755811    &     0.69491650
    \end{matrix}\right].
\end{align*}

\end{exmp}

\begin{exmp}\label{exam-GLQ}
Consider the following three-dimensional state equation
\begin{equation}\label{state-GLQ-exam}
  \left\{
 \begin{aligned}
   \dot{X}(t)&=A\left(\alpha_{t}\right)X(t)+B_{1}\left(\alpha_{t}\right)u_{1}(t)+B_{2}\left(\alpha_{t}\right)u_{2}(t),\\
   X(0)&=x,\quad \alpha(0)=i,
   \end{aligned}
  \right.
\end{equation}
with the cost functional, for any $k=1,2$,
{\scriptsize
\begin{equation}\label{cost-GLQ-exam}
\begin{aligned}
    &J_{k}(x,i;u)\\
    &\!\triangleq \!\mathbb{E}\int_{0}^{\infty}\!
    \left<\!
    \left(\!
    \begin{matrix}
    Q^{k}(\alpha)&S_{1}^{k}(\alpha)^{\top} & S_{2}^{k}(\alpha)^{\top}\\
    S_{1}^{k}(\alpha)&R_{11}^{k}(\alpha) & R_{12}^{k}(\alpha)\\
    S_{2}^{k}(\alpha)&R_{21}^{k}(\alpha) & R_{22}^{k}(\alpha)
    \end{matrix}
    \!\right)\!
    \left(\!
    \begin{matrix}
    X\\
    u_{1}\\
    u_{2}
    \end{matrix}
    \!\right)\!,
    \left(\!
    \begin{matrix}
    X\\
    u_{1}\\
    u_{2}
    \end{matrix}
    \!\right)
    \!\right>dt.
  \end{aligned}
\end{equation}
}
Let
\begin{align*}
 &A(1)\!=\!\left[\begin{matrix}
     -3   &   0   &   0\\
    -1   &   -1   &   0\\
    0   &   0   &   -2
\end{matrix}\right]\!, \quad
A(2)\!=\!\left[\begin{matrix}
     -3   &   0   &   -1\\
    1   &   -2   &   0\\
    0   &   0   &   -1
\end{matrix}\right]\!,\\
&A(3)\!=\!\left[\begin{matrix}
     -2   &   0   &   0\\
    0   &   -1   &   0\\
    0   &   0   &   -1
\end{matrix}\right]\!,\\
&B_{1}(1)\!=\!\left[\begin{matrix}
     2   &   1   \\
    0   &   -2   \\
    0   &   3
\end{matrix}\right]\!,\text{ }
B_{1}(2)\!=\!\left[\begin{matrix}
     2   &   -1   \\
    0   &  1   \\
    2   &   0
\end{matrix}\right]\!,\text{ }
B_{1}(3)\!=\!\left[\begin{matrix}
     1   &   2  \\
    0   &   -1   \\
    1   &   0
\end{matrix}\right]\!,\\
&B_{2}(1)\!=\!\left[\begin{matrix}
     2   &   3   \\
    1   &   -2   \\
    1   &   3
\end{matrix}\right]\!,\text{ }
B_{2}(2)\!=\!\left[\begin{matrix}
     2   &   -1   \\
    0   &   3   \\
    2   &   0
\end{matrix}\right]\!,\text{ }
B_{2}(3)\!=\!\left[\begin{matrix}
     0   &   2  \\
    0   &   4   \\
    1   &   2
\end{matrix}\right]\!,\\
&Q^{1}(1)\!=\!\left[\!\begin{matrix}
 1.79   &   -0.08   &   0.06 \\
-0.08   &   1.69   &   -0.03\\
0.06   &   -0.03   &   1.73
\end{matrix}\!\right]\!,
R_{11}^{1}(1)\!=\!\left[\!\begin{matrix}
   1.55   &   -0.09 \\
   -0.09   &   1.74
\end{matrix}\!\right]\!,\\
&R_{12}^{1}(1)\!=\!\left[\begin{matrix}
 0.15   &   0.03\\
 0.00   &   0.14
\end{matrix}\right]\!,\quad
R_{22}^{1}(1)\!=\!\left[\begin{matrix}
  0.12   &   0.23\\
  0.23   &   0.01
\end{matrix}\right]\!,\\
&S_{1}^{1}(1)^{\top}\!=\!\left[\begin{matrix}
-0.05   &   -0.02   \\
 0.06   &   -0.08 \\
 -0.04   &   -0.01
\end{matrix}\right]\!,\quad
S_{2}^{1}(1)^{\top}\!=\!\left[\begin{matrix}
 0.11   &   0.13\\
  0.14   &   0.10\\
  0.21   &   0.07
\end{matrix}\right]\!,\\
&Q^{1}(2)\!=\!\left[\!\begin{matrix}
1.67   &   -0.05   &   0.00 \\
 -0.05   &   1.80   &   0.12  \\
 0.00   &   0.12   &   1.70
\end{matrix}\!\right]\!,
R_{11}^{1}(2)\!=\!\left[\!\begin{matrix}
   1.58   &   0.04  \\
   0.04   &   1.76
\end{matrix}\!\right]\!,\\
&R_{12}^{1}(2)\!=\!\left[\begin{matrix}
  0.11   &   0.02\\
  0.03   &   0.01
\end{matrix}\right]\!,\quad
R_{22}^{1}(2)\!=\!\left[\begin{matrix}
   0.13   &   0.01\\
  0.01   &   0.02
\end{matrix}\right]\!,\\
&S_{1}^{1}(2)^{\top}\!=\!\left[\begin{matrix}
-0.08   &   0.01\\
-0.04   &  -0.07\\
 0.05   &   0.04
\end{matrix}\right]\!,\quad
S_{2}^{1}(2)^{\top}\!=\!\left[\begin{matrix}
0.23   &   0.12\\
0.13   &   0.01\\
0.24   &   0.01
\end{matrix}\right]\!,\\
&Q^{1}(3)\!=\!\left[\!\begin{matrix}
1.74   &   0.05   &   0.05\\
0.05   &   1.66   &   0.01\\
0.05   &   0.01   &   1.73
\end{matrix}\!\right]\!,
R_{11}^{1}(3)\!=\!\left[\!\begin{matrix}
1.62   &   0.09\\
0.09   &   1.76
\end{matrix}\!\right]\!,\\
&R_{12}^{1}(3)\!=\!\left[\begin{matrix}
0.21   &   0.05\\
0.06   &   0.07
\end{matrix}\right]\!,\quad
R_{22}^{1}(3)\!=\!\left[\begin{matrix}
0.15   &   0.02\\
0.02   &   0.11
\end{matrix}\right]\!,\\
&S_{1}^{1}(3)^{\top}\!=\!\left[\begin{matrix}
-0.09   &   -0.02\\
0.09   &   0.11\\
-0.01   &   -0.03
\end{matrix}\right]\!,\quad
S_{2}^{1}(3)^{\top}\!=\!\left[\begin{matrix}
0.13   &   1.14\\
0.25   &   0.11\\
0.16   &   0.09
\end{matrix}\right]\!,\\
&Q^{2}(1)\!=\!\left[\!\begin{matrix}
1.71   &   0.02   &   -0.02\\
0.02   &   1.61   &   0.03\\
-0.02   &   0.03   &   1.72
\end{matrix}\!\right]\!,
R_{11}^{2}(1)\!=\!\left[\!\begin{matrix}
1.11   &   0.12\\
0.12   &   0.05
\end{matrix}\!\right]\!,\\
&R_{12}^{2}(1)\!=\!\left[\begin{matrix}
0.21   &   0.05\\
0.17   &   0.32
\end{matrix}\right]\!,\quad
R_{22}^{2}(1)\!=\!\left[\begin{matrix}
1.73   &   -0.07\\
-0.07   &   1.73
\end{matrix}\right]\!,\\
&S_{1}^{2}(1)^{\top}\!=\!\left[\begin{matrix}
0.31   &   0.26\\
0.11   &   0.05\\
0.12   &   0.13
\end{matrix}\right]\!,\quad
S_{2}^{2}(1)^{\top}\!=\!\left[\begin{matrix}
0.13   &   0.08\\
-0.02   &   -0.06\\
-0.11  &   0.01
\end{matrix}\right]\!,\\
&Q^{2}(2)\!=\!\left[\!\begin{matrix}
1.84   &   -0.04   &   0.03\\
-0.04   &   1.65   &   -0.03\\
0.03   &   -0.03   &   1.76
\end{matrix}\!\right]\!,
R_{11}^{2}(2)\!=\!\left[\!\begin{matrix}
0.07   &   0.15\\
0.15   &   0.28
\end{matrix}\!\right]\!,\\
&R_{12}^{2}(2)\!=\!\left[\begin{matrix}
0.21   &   0.07\\
0.71   &   0.15
\end{matrix}\right]\!,\quad
R_{22}^{2}(2)\!=\!\left[\begin{matrix}
1.61   &   0.05\\
0.05   &   1.64
\end{matrix}\right]\!,\\
&S_{1}^{2}(2)^{\top}\!=\!\left[\begin{matrix}
0.15   &   0.18\\
0.16   &   0.17\\
0.11   &   0.31
\end{matrix}\right]\!,\quad
S_{2}^{2}(2)^{\top}\!=\!\left[\begin{matrix}
-0.11   &   0.04\\
0.07   &   0.02\\
0.07   &   0.02
\end{matrix}\right]\!,\\
&Q^{2}(3)\!=\!\left[\!\begin{matrix}
1.59   &  -0.07   &   -0.01\\
-0.07   &   1.77   &   0.02\\
-0.01   &   0.02   &   1.66
\end{matrix}\!\right]\!,
R_{11}^{2}(3)\!=\!\left[\!\begin{matrix}
0.35   &   0.17\\
0.17   &   0.21
\end{matrix}\!\right]\!,\\
&R_{12}^{2}(3)\!=\!\left[\begin{matrix}
0.09   &   0.12\\
0.37   &   0.15
\end{matrix}\right]\!,\quad
R_{22}^{2}(3)\!=\!\left[\begin{matrix}
1.73   &   -0.04\\
-0.04   &   1.75
\end{matrix}\right]\!,\\
&S_{1}^{2}(3)^{\top}\!=\!\left[\begin{matrix}
0.11   &   0.31\\
0.17   &   0.28\\
0.16   &   0.41
\end{matrix}\right]\!,\quad
S_{2}^{2}(3)^{\top}\!=\!\left[\begin{matrix}
-0.04   &   0.02\\
-0.03   &   -0.05\\
0.10   &   -0.07
\end{matrix}\right]\!.
\end{align*}
Solving the corresponding cross-coupled CAREs \eqref{CAREs-GLQ-1}-\eqref{CAREs-GLQ-2}, we obtain
\begin{align*}
     P_{1}(1)&=\left[\begin{matrix}
     0.21584298  & -0.06035326  &  -0.08509888\\
    -0.06035326  &  0.36829891  &   0.07611724\\
    -0.08509888  &  0.07611724  &   0.31182214
    \end{matrix}\right],\\
     P_{1}(2)&=\left[\begin{matrix}
     0.25749695  &  0.03556257  &  -0.14493485\\
     0.03556257  &  0.26164453  &  -0.01161552\\
    -0.14493485  & -0.01161552  &   0.47583037
    \end{matrix}\right],\\
     P_{1}(3)&=\left[\begin{matrix}
     0.28959329  & -0.15275789  &  -0.07306250\\
    -0.15275789  &  0.28516579  &  -0.16017357\\
    -0.07306250  & -0.16017357  &   0.50230499
    \end{matrix}\right],\\
    P_{2}(1)&=\left[\begin{matrix}
     0.21105400  & -0.04533077  &  -0.07477181\\
    -0.04533077  &  0.44894270  &   0.13781734\\
    -0.07477181  &  0.13781734  &   0.30765576
    \end{matrix}\right],\\
     P_{2}(2)&=\left[\begin{matrix}
     0.29559748  &  0.05746950  &  -0.14594113\\
     0.05746950  &  0.30207529  &  -0.04821299\\
    -0.14594113  & -0.04821299  &   0.49999012
    \end{matrix}\right],\\
     P_{2}(3)&=\left[\begin{matrix}
     0.22967392  & -0.00025983  &  -0.06922013\\
    -0.00025983  &  0.37356702  &  -0.08557568\\
    -0.06922013  & -0.08557568  &   0.52651832
    \end{matrix}\right],
\end{align*}
with the corresponding residuals are given by
\begin{align*}
&\|E_{1}(\mathbf{P_{1}})\|\!= \!7.3356\!\times\! 10^{-8}\!,\text{ }
\|E_{2}(\mathbf{P_{1}})\|\!=\! 7.8884\!\times \! 10^{-8}\!,\\
&\|E_{3}(\mathbf{P_{1}})\|\!=\!5.7103\!\times\! 10^{-8}\!,\text{ }
\|E_{1}(\mathbf{P_{2}})\|\!=\!5.6206\!\times\! 10^{-8}\!,\\
&\|E_{2}(\mathbf{P_{2}})\|\!=\!6.3018\!\times\! 10^{-8}\!,\text{ }
\|E_{3}(\mathbf{P_{2}})\|\!=\!1.0050\!\times\! 10^{-7}\!.
\end{align*}
Then, by Theorem \ref{thm-GLQ-closed} and Corollary \ref{coro-GLQ-0-closed}, we obtain the following closed-loop Nash equilibrium point of Problem (M-GLQ)$^{0}$ via solving the \eqref{CAREs-GLQ-constraint} as follows:
\begin{align*}
    \widehat{\Theta}_{1}(1)&=\left[\begin{matrix}
    -0.21818337    &   0.08325908   &    0.12663298\\
    -0.01957014    &   0.35184937   &   -0.37375256
    \end{matrix}\right],\\
     \widehat{\Theta}_{1}(2)&=\left[\begin{matrix}
    -0.08425975   &   0.00388356   &    -0.42138102\\
     0.12484096   &  -0.08578065   &    -0.08234423
    \end{matrix}\right],\\
     \widehat{\Theta}_{1}(3)&=\left[\begin{matrix}
    -0.07049709   &   0.14753455    &   -0.20478132\\
    -0.39994949   &   0.29667540    &    0.04339824
    \end{matrix}\right],\\
    \widehat{\Theta}_{2}(1)&=\left[\begin{matrix}
    -0.23484315    &   -0.30720226   &    -0.09375016\\
    -0.33455232    &    0.31339120   &    -0.18861879
    \end{matrix}\right],\\
     \widehat{\Theta}_{2}(2)&=\left[\begin{matrix}
    -0.16313607    &   -0.00144350   &    -0.39277560\\
     0.04787672    &   -0.52200546   &     0.02450300
    \end{matrix}\right],\\
     \widehat{\Theta}_{2}(3)&=\left[\begin{matrix}
    0.14883148     &   -0.02196658    &    -0.36736956\\
   -0.15169310     &   -0.76324639    &    -0.28509967
    \end{matrix}\right].
\end{align*}
\end{exmp}

\bibliography{references}
\bibliographystyle{AIMS}

\end{document}